\documentclass[12pt,reqno]{amsart}
\usepackage{amsmath}
\usepackage{latexsym,amssymb}
\usepackage{mathrsfs}
\usepackage{enumerate}
\usepackage{bm}
\allowdisplaybreaks[1]

\newtheorem{lemma}{Lemma}[section]
\newtheorem{theorem}[lemma]{Theorem}
\newtheorem{thme}[]{Theorem}

\newtheorem{proposition}[lemma]{Proposition}
\newtheorem{corollary}[lemma]{Corollary}
\newtheorem{question}[lemma]{Question}
\theoremstyle{definition}
\newtheorem{definition}[lemma]{Definition}
\newtheorem{example}[lemma]{Example}
\newtheorem{remark}[lemma]{Remark}

\numberwithin{equation}{section}

\newcommand{\bdf}{\begin{definition}}
\newcommand{\edf}{\end{definition}}
\newcommand{\blem}{\begin{lemma}}
\newcommand{\elem}{\end{lemma}}
\newcommand{\bthm}{\begin{theorem}}
\newcommand{\ethm}{\end{theorem}}
\newcommand{\bpf}{\begin{proof}}
\newcommand{\epf}{\end{proof}}
\newcommand{\bprop}{\begin{proposition}}
\newcommand{\eprop}{\end{proposition}}
\newcommand{\bcor}{\begin{corollary}}
\newcommand{\ecor}{\end{corollary}}
\newcommand{\brem}{\begin{remark}}
\newcommand{\erem}{\end{remark}}
\newcommand{\bquest}{\begin{question}}
\newcommand{\equest}{\end{question}}
\newcommand{\bex}{\begin{example}}
\newcommand{\eex}{\end{example}}

\newcommand{\benu}{\begin{enumerate}\renewcommand{\labelenumi}{{\rm (\arabic{enumi})}}\renewcommand{\itemsep}{0pt}}
\newcommand{\eenu}{\end{enumerate}}

\newcommand{\M}{\mathbb{M}}
\newcommand{\N}{\mathbb{N}}

\newcommand{\Q}{\mathbb{Q}}
\newcommand{\R}{\mathbb{R}}
\newcommand{\C}{\mathbb{C}}

\newcommand{\bB}{\mathbb{B}}

\def\cA{\mathcal{A}}

\def\cF{\mathcal{F}}
\def\cG{\mathcal{G}}

\DeclareMathOperator{\Aut}{Aut}

\DeclareMathOperator{\tr}{tr}
\DeclareMathOperator{\Tr}{Tr}

\newcommand{\stand}{(M, H, J, P)}
\newcommand{\mat}{\mathbb{M}_n}
\newcommand{\e}{\varepsilon}

\def\al{\alpha}
\def\hal{{\hat{\al}}}

\def\be{\beta}
\def\ga{\gamma}
\def\de{\delta}

\def\la{\lambda}

\def\vep{\varepsilon}

\def\ps{{\psi}}

\def\vph{\varphi}

\def\om{\omega}
\def\si{\sigma}

\def\th{\theta}

\def\De{\Delta}

\def\La{\Lambda}
\def\Ph{\Phi}
\def\Ps{\Psi}

\def\subs{\subset}
\def\ovl{\overline}
\def\oti{\otimes}

\DeclareMathOperator{\id}{id}

\begin{document}
\title{Haagerup approximation property and positive cones associated with a von Neumann algebra}

\author[R. Okayasu]{Rui Okayasu$^1$}
\address{$^1$
Department of Mathematics Education, Osaka Kyoiku University,
Osaka \mbox{582-8582},
JAPAN}
\email{rui@cc.osaka-kyoiku.ac.jp}

\author[R. Tomatsu]{Reiji Tomatsu$^2$}
\address{$^2$
Department of Mathematics, Hokkaido University,
Hokkaido \mbox{060-0810},
JAPAN}
\email{tomatsu@math.sci.hokudai.ac.jp}

\date{\today}
\thanks{The first author was partially supported by JSPS KAKENHI Grant Number 25800065. The second author was partially supported by JSPS KAKENHI Grant Number 24740095.}

\maketitle

\begin{abstract}
We introduce the notion of the $\al$-Haagerup
approximation property for $\al\in[0,1/2]$
using a one-parameter family of positive cones
studied by Araki
and show that the $\al$-Haagerup approximation
property actually does not depend
on a choice of $\al$.
This enables us to prove the fact that
the Haagerup approximation properties introduced in two ways
are actually equivalent,
one in terms of the standard form
and the other in terms of completely positive maps.
We also discuss the $L^p$-Haagerup approximation property
for a non-commutative $L^p$-spaces associated with a von Neumann algebra
($1<p<\infty$)
and show the independence of the $L^p$-Haagerup approximation property on $p$.
\end{abstract}

%%%%%%%%%%%%%%%%%%%%%%%%%%%%%%%%%%%%%%%%%%%%%%%%%%%%%%%%%%%%%%%%%%%%%%%%%%%%%%%%%%%%%%%%%%%%%%%%%%%%%%%%%%%%%%%%%%%%%%%%%%%%%%%%%%%%%%%%%%%%%%%%%%%%%%%%%%%%%%%%%%%%%%%%%%%%%%%%%%%%%%%%%%%%%%%%%%%%%%%%%%%%%%%%%%%%%%%%%%%%%%%%%%%%%%%%%%%%%

\section{Introduction}

%%%%%%%%%%%%%%%%%%%%%%%%%%%%%%%%%%%%%%%%%%%%%%%%%%%%%%%%%%%%%%%%%%%%%%%%%%%%%%%%%%%%%%%%%%%%%%%%%%%%%%%%%%%%%%%%%%%%%%%%%%%%%%%%%%%%%%%%%%%%%%%%%%%%%%%%%%%%%%%%%%%%%%%%%%%%%%%%%%%%%%%%%%%%%%%%%%%%%%%%%%%%%%%%%%%%%%%%%%%%%%%%%%%%%%%%%%%%%

This is a continuation of our previous work \cite{ot}
on the Haagerup approximation property (HAP)
for a von Neumann algebra.
The origin of the HAP
is
%in
the remarkable paper \cite{haa3},
where U. Haagerup proved that 
the reduced group $\mathrm{C}^*$-algebra of the non-amenable free group
has Grothendieck's metric approximation property.
After his work,
M. Choda \cite{cho} showed that 
a discrete group
has the HAP
%the Haagerup approximation property
if and only if
its group von Neumann algebra has 
a certain von Neumann algebraic approximation property 
with respect to the natural faithful normal tracial state.
Furthermore, P. Jolissaint \cite{jol} studied 
the HAP
%the Haagerup approximation property
in the framework of finite von Neumann algebras.
In particular, it was proved that
it does not depend on the choice of a faithful normal tracial state. 

In the last few years,
the Haagerup type approximation property for quantum groups
with respect to the Haar states 
was actively investigated by many researchers 
(e.g.\ \cite{br1,br2,dfsw,cfy,kv,le}). 
The point here is that 
the Haar state on a quantum group
is not necessarily tracial,
and so to fully understand the HAP
%the Haagerup approximation property
for quantum groups,
we need to characterize this property
in the framework of arbitrary von Neumann algebras.

In the former work \cite{ot},
we introduce the notion of
the HAP
%the Haagerup approximation property
for arbitrary von Neumann algebras
in terms of the standard form.
Namely,
the HAP means the existence of 
contractive completely positive compact operators
on the standard Hilbert space
which are approximating to the identity.
In \cite{cs},
M. Caspers and A. Skalski 
independently introduce 
the notion 
of the HAP
based on the existence of 
completely positive maps approximating to the identity
with respect to a given faithful normal semifinite weight
such that the associated implementing operators
on the GNS Hilbert space are compact.

Now one may wonder
whether these two approaches are different or not.
Actually, by combining several results in \cite{ot} and \cite{cs},
it is possible to show that
these two formulations are equivalent.
(See \cite{cost}, \cite[Remark 5.8]{ot} for details.)
This proof, however, relies on
the permanence results of the HAP
for a core von Neumann algebra.
%and is unsatisfactory.
One of our purposes in the present paper
is to give a simple and direct proof
for the above mentioned question.

Our strategy is to use the positive cones
due to H. Araki.
%Following the development of the Tomita--Takesaki theory,
He introduced in \cite{ara}
a one-parameter family of positive cones $P^\alpha$
with a parameter $\alpha$ in the interval $[0, 1/2]$
that is associated with a von Neumann algebra 
admitting a cyclic and separating vector.
This family is ``interpolating''
the three distinguished cones 
$P^0$, $P^{1/4}$ and $P^{1/2}$,
which are also denoted by 
$P^\sharp$, $P^\natural$ and $P^\flat$
in the literature \cite{t2}.
Among them, the positive cone
$P^\natural$ at the middle point
plays remarkable roles
in the theory of the standard representation \cite{ara,co1,haa1}.
See \cite{ara,ko1,ko2} for comprehensive
studies of that family.

In view of the positive cones $P^\al$,
on the one hand, 
our definition of the HAP
is, of course, related with $P^\natural$.  
On the other hand,
the associated $L^2$-GNS implementing operators
in the definition due to Caspers and Skalski
are, in fact, ``completely positive''
with respect to $P^\sharp$.
Motivated by these facts,
we will introduce the notion of the ``interpolated'' 
HAP called $\al$-HAP
and prove the following result
(Theorem \ref{thm:equiv}):

\begin{thme}
A von Neumann algebra $M$ 
has the $\alpha$-HAP 
for some $\alpha\in[0, 1/2]$
if and only if 
$M$ has the $\alpha$-HAP 
for all $\alpha\in[0, 1/2]$
\end{thme}

As a consequence, 
it gives a direct proof that
two definitions of the HAP introduced in \cite{cs,ot}
are equivalent.

In the second part of the present paper,
we discuss the Haagerup approximation property
for non-commutative $L^p$-spaces ($1<p<\infty$)
\cite{am,haa2,han,izu,ko3,te1,te2}.
One can introduce the natural notion of
the complete positivity
of operators on $L^p(M)$,
and hence we will define the HAP
called the $L^p$-HAP
when there exists a net of
completely positive compact operators
approximating to the identity on $L^p(M)$.
Since $L^2(M)$ is the standard form of $M$,
it follows from the definition that
a von Neumann algebra $M$ has the HAP
if and only if $M$ has the $L^2$-HAP.
Furthermore,
by using the complex interpolation method
due to A. P. Calder\'{o}n \cite{ca},
we can show the following result
(Theorem \ref{thm:L^p-HAP3}):

\begin{thme}
Let $M$ be a von Neumann algebra.
Then the following statements are equivalent:
\benu
\item $M$ has the HAP;
\item $M$ has the $L^p$-HAP for all $1<p<\infty$;
\item $M$ has the $L^p$-HAP for some $1<p<\infty$.
\eenu
\end{thme}

We remark that a von Neumann algebra $M$ 
has the completely positive approximation property (CPAP)
if and only if $L^p(M)$ has the CPAP for some/all $1\leq p<\infty$.
In the case where $p=1$, this is proved by E. G. Effros and E. C. Lance in \cite{el}.
In general, this is due to M. Junge, Z-J. Ruan and Q. Xu in \cite{jrx}.
Therefore Theorem B is the HAP version of this result.

%%%%%%%%%%%%%%%%%%%%%%%%%%%%%%%%%%%%%%%%%%%%%%%%%%%%%%%%%%%%%%%%%%%%%%%%%%%%%%%%%%%%%%%%%%%%%%%%%%%%%%%%%%%%%%%%%%%%%%%%%%%%%%%%%%%%%%%%%%%%%%%%%%%%%%%%%%%%%%%%%%%%%%%%%%%%%%%%%%%%%%%%%%%%%%%%%%%%%%%%%%%%%%%%%%%%%%%%%%%%%%%%%%%%%%%%%%%%

\vspace{10pt}
\noindent
\textbf{ Acknowledgments.}
The authors would like to thank
Marie Choda and Yoshikazu Katayama for their encouragement
and fruitful discussion,
and Martijn Caspers and Adam Skalski
for valuable comments on our work.
They also would like to thank
Yoshimichi Ueda for stimulating discussion.

%%%%%%%%%%%%%%%%%%%%%%%%%%%%%%%%%%%%%%%%%%%%%%%%%%%%%%%%%%%%%%%%%%%%%%%%%%%%%%%%%%%%%%%%%%%%%%%%%%%%%%%%%%%%%%%%%%%%%%%%%%%%%%%%%%%%%%%%%%%%%%%%%%%%%%%%%%%%%%%%%%%%%%%%%%%%%%%%%%%%%%%%%%%%%%%%%%%%%%%%%%%%%%%%%%%%%%%%%%%%%%%%%%%%%%%%%%%%%

\section{Preliminaries}

%%%%%%%%%%%%%%%%%%%%%%%%%%%%%%%%%%%%%%%%%%%%%%%%%%%%%%%%%%%%%%%%%%%%%%%%%%%%%%%%%%%%%%%%%%%%%%%%%%%%%%%%%%%%%%%%%%%%%%%%%%%%%%%%%%%%%%%%%%%%%%%%%%%%%%%%%%%%%%%%%%%%%%%%%%%%%%%%%%%%%%%%%%%%%%%%%%%%%%%%%%%%%%%%%%%%%%%%%%%%%%%%%%%%%%%%%%%%

We first fix the notation 
and recall several facts studied in \cite{ot}.
Let $M$ be a von Neumann algebra. 
We denote by $M_{\mathrm{sa}}$ and $M^+$, 
the set of all self-adjoint elements 
and all positive elements in $M$, respectively. 
We also denote by $M_*$ and $M_*^+$, 
the space of all normal linear functionals 
and all positive normal linear functionals on $M$, respectively.
The set of faithful normal semifinite (f.n.s.) weights
is denoted by $W(M)$.
Recall the definition of a standard form
of a von Neumann algebra.

\bdf[{\cite[Definition 2.1]{haa1}}]

Let $\stand$ be a quadruple,
where $M$ denotes a von Neumann algebra, 
$H$ a Hilbert space on which $M$ acts, 
$J$ a conjugate-linear isometry on $H$ with $J^2=1_H$, 
and $P\subset H$
a closed convex cone which is self-dual,
i.e., $P=P^\circ$, where
$P^\circ
:=\{\xi\in H \mid \langle\xi, \eta\rangle\geq 0
\ \text{for } \eta\in H\}$.
Then $\stand$ is called a {\em standard form} 
if the following conditions are satisfied: 

\benu
\item
$J M J= M'$;

\item
$J\xi=\xi$ for any $\xi\in P$;

\item
$aJaJ P\subset P$ for any $a\in M$;

\item
$JcJ=c^*$ for any $c\in \mathcal{Z}(M):= M\cap M'$.

\eenu

\edf

\brem

In \cite{ah1}, Ando and Haagerup proved 
that the condition (4) in the above definition can be removed. 

\erem

We next introduce that each f.n.s.\ weight $\vph$ 
gives a standard form. 
We refer readers
to the book of Takesaki \cite{t2}
for details. 
Let $M$ be a von Neumann algebra with $\vph\in W(M)$.
We write 
\[
n_\varphi
:=\{x\in M \mid \varphi(x^*x)<\infty\}.
\]
Then $H_\varphi$ is the completion of $n_\varphi$ with respect to the norm 
\[
\|x\|_\varphi^2
:=\varphi(x^*x)\quad \text{for } x\in n_\vph.
\]
We write the canonical injection 
$\Lambda_\varphi\colon n_\varphi\to H_\varphi$. 

Then 
\[
\mathcal{A}_\varphi
:=\Lambda_\varphi(n_\varphi\cap n_\varphi^*)
\] 
is an achieved left Hilbert algebra with the multiplication 
\[
\Lambda_\varphi(x)\cdot\Lambda_\varphi(x)
:=\Lambda_\varphi(xy)\quad \text{for } x\in n_\varphi\cap n_\varphi^*
\]
and the involution 
\[
\Lambda_\varphi(x)^\sharp
:=\Lambda_\varphi(x^*)\quad \text{for } x\in n_\varphi\cap n_\varphi^*.
\]
Let $\pi_\varphi$ be the corresponding representation of $M$ on $H_\varphi$. 
We always identify $M$ with $\pi_\varphi(M)$. 

We denote by $S_\varphi$ the closure 
of the conjugate-linear operator $\xi\mapsto\xi^\sharp$ 
on $H_\varphi$, 
which has the polar decomposition 
\[
S_\varphi
=J_\varphi\Delta_\varphi^{1/2},
\]
where $J_\varphi$ is the modular conjugation 
and $\Delta_\varphi$ is the modular operator. 
The modular automorphism group $(\si_t^\vph)_{t\in\R}$
is given by
\[
\si_t^\vph(x):=\De_\vph^{it}x\De_\vph^{-it}
\quad\text{for } x\in M.
\]
For $\vph\in W(M)$, 
we denote the centralizer of $\vph$
by 
\[
M_\vph:=\{x\in M \mid \si_t^\vph(x)=x\ \text{for } t\in\R\}.
\]

Then we have a self-dual positive cone
\[
P_\varphi^\natural
:=\overline{\{\xi(J_\varphi\xi) \mid \xi\in\mathcal{A}_\varphi\}}
\subset H_\varphi.
\]
Note that $P_\vph^\natural$ is given by the closure of the set
of $\La_\vph(x\si_{i/2}^\vph(x)^*)$,
where $x\in \cA_\vph$ is entire
with respect to $\si^\vph$.

Therefore
the quadruple $(M, H_\varphi, J_\varphi, P_\varphi^\natural)$
is a standard form. 
Thanks to \cite[Theorem 2.3]{haa1},
a standard form is, in fact, unique
up to a spatial isomorphism,
and so it is independent to the choice
of an f.n.s.\ weight $\varphi$.

Let us consider the $n\times n$ matrix algebra $\mat$ 
and the normalized trace $\mathrm{tr}_n$. 
The algebra $\mat$ becomes a Hilbert space
with the inner product
$\langle x, y\rangle:=\mathrm{tr}_n(y^*x)$
for $x,y\in\mat$.
We write the canonical involution 
$J_{\mathrm{tr}_n}\colon x\mapsto x^*$ for $x\in \mat$. 
Then the quadruple $(\mat, \mat, J_{\mathrm{tr}_n}, \mat^+)$ 
is a standard form.
In the following,
for a Hilbert space $H$,
$\mat(H)$ denotes the tensor product Hilbert space
$H\oti \mat$.

\bdf[{\cite[Definition 2.2]{mt}}]
Let $\stand$ be a standard form and $n\in\N$. 
A matrix $[\xi_{i, j}]\in\mat(H)$ is said to be {\em positive} if 
\[
\sum_{i, j=1}^nx_iJx_jJ\xi_{i, j}\in P
\quad
\text{for all }
x_1,\dots,x_n\in M.
\]
We denote by $P^{(n)}$ the set of all positive matrices $[\xi_{i, j}]$ in $\mat(H)$. 
\edf

\bprop[{\cite[Proposition 2.4]{mt}}, {\cite[Lemma 1.1]{sw1}}]
Let $\stand$ be a standard form and $n\in\N$. 
Then $(\mat(M), \mat(H), J\otimes J_{\mathrm{tr}_n}, P^{(n)})$ 
is a standard form.
\eprop

Next, we will introduce the complete positivity
of a bounded operator between standard Hilbert spaces.

\bdf
\label{defn:cpop}
Let $(M_1, H_1, J_1, P_1)$ and $(M_2, H_2, J_2, P_2)$ be two standard forms. 
We will say that a bounded linear (or conjugate-linear) operator 
$T\colon H_1\to H_2$ is {\em completely positive} 
if $(T\otimes1_{\mat})P_1^{(n)}\subset P_2^{(n)}$
for all $n\in\N$.
\edf

\bdf[{\cite[Definition 2.7]{ot}}]
\label{defn:HAP}
A W$^*$-algebra $M$
has the {\em Haagerup approximation property} (HAP)
if there exists a standard form $(M, H, J, P)$
and a net of contractive completely positive (c.c.p.)
compact operators $T_n$ on $H$ 
such that $T_n\to 1_{H}$ in the strong topology.
\edf

Thanks to \cite[Theorem 2.3]{haa1},
this definition does not depend on the choice of a standard from.
We also remark that the weak convergence
of a net $T_n$ in the above definition is sufficient.
In fact, we can arrange a net $T_n$ such that 
$T_n\to 1_H$ in the strong topology
by taking suitable convex combinations.

In the case where $M$ is $\sigma$-finite 
with a faithful state $\varphi\in M_*^+$.
We denote by $(H_\vph, \xi_\vph)$ 
the GNS Hilbert space 
with the cyclic and separating vector
associated with $(M, \vph)$.
If $M$ has the HAP,
then we can recover a net of c.c.p.\ maps on $M$
approximating to the identity
with respect to $\vph$
such that the associated implementing operators
on $H_\vph$ are compact.

\bthm[{\cite[Theorem 4.8]{ot}}]
\label{thm:sigma-finite}
Let $M$ be a $\sigma$-finite von Neumann algebra 
with a faithful state $\varphi\in M_*^+$. 
Then $M$ has the HAP if and only if 
there exists a net of normal c.c.p.\ maps $\Phi_n$ on $M$ such that
\begin{itemize}
\item
$\varphi\circ\Phi\leq \varphi$;

\item
$\Phi_n\to\mathrm{id}_{M}$ in the point-ultraweak topology;

\item
The operator defined below
is c.c.p.\ compact on $H_\varphi$
and $T_n\to 1_{H_\varphi}$ in the strong topology:
\[
T_n(\Delta_\varphi^{1/4}x\xi_\varphi)
=\Delta_\varphi^{1/4}\Phi_n(x)\xi_\varphi
\ \text{for } x\in M.
\] 
\end{itemize}
\ethm

This translation of the HAP looks similar to the following HAP
introduced by Caspers and Skalski in \cite{cs}.

\bdf[{\cite[Definition 3.1]{cs}}]
\label{defn:CS1}
Let $M$ be a von Neumann algebra
with $\vph\in W(M)$.
We will say that $M$ has 
the \emph{Haagerup approximation property
with respect to $\vph$}
in the sense of \cite{cs}
(CS-HAP$_\vph$)
if
there exists a net of normal c.p.\ maps $\Ph_n$ on $M$
such that
\begin{itemize}
\item
$\vph\circ\Ph_n\leq\vph$;

%\item
%$\Ph_n\to\id_M$ in the point-ultraweak topology;

\item
The operator $T_n$ defined below is compact
and $T_n\to 1_{H_\varphi}$ in the strong topology:
\[
T_n\La_\vph(x):=\La_\vph(\Ph_n(x))
\quad
\mbox{for }
x\in n_\vph.
\]
\end{itemize}
\edf

Here are two apparent differences
between Theorem \ref{thm:sigma-finite}
and Definition \ref{defn:CS1},
that is,
the existence of $\De_\vph^{1/4}$ of course,
and the assumption on the contractivity of $\Ph_n$'s.
Actually,
it is possible to show that 
the notion of the CS-HAP$_\vph$ does not depend 
on the choice of $\vph$ \cite[Theorem 4.3]{cs}.
Furthermore we can take contractive $\Ph_n$'s. 
(See Theorem \ref{thm:cpmap}.)
The proof of the weight-independence presented in \cite{cs}
relies on a crossed product work.
Here, let us present a direct proof
of the weight-independence of the CS-HAP.

\blem[{\cite[Theorem 4.3]{cs}}]
\label{lem:CS-indep}
The CS-HAP is the weight-free property.
Namely,
let $\vph, \ps\in W(M)$.
Then $M$ has the CS-HAP$_\vph$
if and only if
$M$ has the CS-HAP$_\ps$.
\elem

\bpf
%Let $\vph_0$ be the model weight given
%in the proof of Lemma \ref{lem:if}.
%Then for each $e_i\in M_{\vph_0}$,
%we can construct a net of c.c.p.\ maps
%$\Ph_n^i$ satisfying the required property.
%Then their summation does the job,
%that is,
%the statement for $\vph_0$ holds.
%
Suppose that
$M$ has the CS-HAP$_\vph$.
Let $\Ph_n$ and $T_n$ be as in the statement
of Definition \ref{defn:CS1}.
%Let us take an arbitrary f.n.s.\ weight $\ps$ on $M$.
Note that an arbitrary $\ps\in W(M)$ is obtained from $\vph$
by combining the following four operations:
\begin{enumerate}
\item
$\vph\mapsto \vph\oti\Tr$,
where $\Tr$ denotes the canonical tracial weight
on $\bB(\ell^2)$;

\item
$\vph\mapsto \vph_e$,
where $e\in M_\vph$ is a projection;

\item
$\vph\mapsto \vph\circ\al$,
$\al\in \Aut(M)$;

\item
$\vph\mapsto \vph_h$,
where $h$ is a non-singular positive operator
affiliated with $M_\vph$
and $\vph_h(x):=\vph(h^{1/2}xh^{1/2})$ for $x\in M^+$.
\end{enumerate}
For its proof,
see the proof of \cite[Th\'{e}or\`{e}me 1.2.3]{co1}
or \cite[Corollary 5.8]{st}.
Hence it suffices to consider each operation.

(1)
Let $\ps:=\vph\oti\Tr$.
Take an increasing net of finite rank projections $p_n$
on $\ell^2$.
Then $\Ph_n\oti (p_n\cdot p_n)$ does the job,
where $p_n\cdot p_n$ means the map $x\mapsto p_n xp_n$.

(2)
Let $e\in M_\vph$ be a projection.
Set $\ps:=\vph_e$
and $\Ps_n:=e\Ph_n(e\cdot e)e$.
Then we have $\ps\circ\Ps_n\leq\ps$.
Indeed, for $x\in (eMe)_+$,
we obtain
\[
\ps(x)
=
\vph(exe)
\geq
\vph(\Ph_n(exe))
\geq
\vph(e\Ph_n(exe)e)
=
\ps(\Ps_n(x)).
\]
Moreover for $x\in n_\vph$,
we have
\begin{align*}
\La_{\vph_e}(\Ps_n(exe))
&=
eJeJ \La_\vph(\Ph_n(exe))
\\
&=
eJeJ T_n
\La_\vph(exe)
\\
&=
eJeJ T_neJeJ
\La_{\vph_e}(exe).
\end{align*}
Since $eJeJ T_n eJeJ$
is compact, we are done.

(3)
Let $\ps:=\vph\circ\al$.
Regard as $H_\psi=H_\vph$ by putting $\La_\ps=\La_\vph\circ\al$.
Then we obtain the canonical unitary implementation $U_\al$
which maps $\La_\vph(x)\mapsto\La_\ps(\al^{-1}(x))$
for $x\in n_\vph$.
Set $\Ps_n:=\al^{-1}\circ\Phi_n\circ\al$.
Then we have
\[
\ps(x)=\vph(\al(x))
\geq\vph(\Phi_n(\al(x)))
=\ps(\Ps_n(x))
\quad\text{for } x\in M^+, 
\]
and 
\[
U_\al T_nU_\al^*\La_\ps(x)
=U_\al T_n\La_\vph(\al(x))
=U_\al \La_\vph(\Phi_n(\al(x)))
=\La_\ps(\Ps(x))
\quad\text{for } x\in n_\vph.
\]
Since $U_\al T_nU_\al^*$ is compact, we are done.

(4)
This case is proved in \cite[Proposition 4.2]{cs}.
Let us sketch out its proof for readers' convenience.
Let $e(\cdot)$ be the spectral resolution of $h$
and put $e_n:=e([1/n,n])$ for $n\in\N$.
Considering $\vph_{he_n}$,
we may and do assume that
$h$ is bounded and invertible
by \cite[Lemma 4.1]{cs}.
Put
$\Ps_n(x):=h^{-1/2}\Ph_n(h^{1/2}xh^{1/2})h^{-1/2}$
for $x\in M$.
Then we have
$\vph_h\circ\Ps_n\leq\vph_h$,
and
the associated implementing operator
is given by
$h^{-1/2}T_nh^{1/2}$, which is compact.
\epf

%%%%%%%%%%%%%%%%%%%%%%%%%%%%%%%%%%%%%%%%%%%%%%%%%%%%%%%%%%%%%%%%%%%%%%%%%%%%%%%%%%%%%%%%%%%%%%%%%%%%%%%%%%%%%%%%%%%%%%%%%%%%%%%%%%%%%%%%%%%%%%%%%%%%%%%%%%%%%%%%%%%%%%%%%%%%%%%%%%%%%%%%%%%%%%%%%%%%%%%%%%%%%%%%%%%%%%%%%%%%%%%%%%%%%%%%%%%%%

\section{Haagerup approximation property and positive cones}

%%%%%%%%%%%%%%%%%%%%%%%%%%%%%%%%%%%%%%%%%%%%%%%%%%%%%%%%%%%%%%%%%%%%%%%%%%%%%%%%%%%%%%%%%%%%%%%%%%%%%%%%%%%%%%%%%%%%%%%%%%%%%%%%%%%%%%%%%%%%%%%%%%%%%%%%%%%%%%%%%%%%%%%%%%%%%%%%%%%%%%%%%%%%%%%%%%%%%%%%%%%%%%%%%%%%%%%%%%%%%%%%%%%%%%%%%%%%%

In this section, we generalize
the HAP using
a one-parameter family of positive cones 
parametrized by $\alpha\in [0, 1/2]$, 
which is introduced by Araki in \cite{ara}.
Let $M$ be a von Neumann algebra
and $\vph\in W(M)$.

%%%%%%%%%%%%%%%%%%%%%%%%%%%%%%%%%%%%%%%%%%%%%%%%%%%%%%%%%%%%%%%%%%%%%%%%%%%%%%%

\subsection{Complete positivity associated with positive cones}

%%%%%%%%%%%%%%%%%%%%%%%%%%%%%%%%%%%%%%%%%%%%%%%%%%%%%%%%%%%%%%%%%%%%%%%%%%%%%%%

Recall that $\cA_\vph$ is the associated left Hilbert algebra.
Let us consider the following positive cones:
\[
P_\varphi^\sharp
:=\overline{\{ \xi\xi^\sharp \mid \xi\in \cA_\varphi\}},
\quad
P_\varphi^\natural
:=\overline{\{\xi(J_\varphi\xi) \mid \xi\in\cA_\varphi\}},
\quad
P_\varphi^\flat:=\overline{\{ \eta\eta^\flat \mid \xi\in \cA_\varphi'\}}
\]
Then $P_\vph^\sharp$ is contained in $D(\De_\vph^{1/2})$,
the domain of $\De_\vph^{1/2}$.

\bdf[cf. {\cite[Section 4]{ara}}]
For $\alpha\in[0, 1/2]$,
we will define the positive cone $P_\varphi^\alpha$
by the closure of $\Delta_\varphi^\alpha P_\varphi^\sharp$.
\edf

Then $P_\vph^\al$ has
the same properties as in \cite[Theorem 3]{ara}:

\benu
\item $P_\varphi^\alpha$ is the closed convex cone
invariant under $\Delta_\varphi^{it}$;

\item $P_\varphi^\alpha
\subset D(\Delta_\varphi^{1/2-2\alpha})$ 
and $J_\varphi\xi=\Delta_\varphi^{1/2-2\alpha}\xi$ 
for $\xi\in P_\varphi^\alpha$;

\item $J_\varphi P_\varphi^\alpha=P_\varphi^{\hat{\al}}$,
where $\hat{\al}:=1/2-\alpha$;

\item $P_\varphi^{\hat{\al}}
=\{\eta\in H_\varphi \mid 
\langle \eta, \xi\rangle\geq 0 
\ \text{for }
\xi\in P_\varphi^\alpha\}$;

\item $P_\varphi^\alpha=\Delta_\varphi^{\alpha-1/4}(P_\varphi^{1/4}\cap D(\Delta_\varphi^{\alpha-1/4}))$;
\item $P_\varphi^\natural=P_\varphi^{1/4}$ and $P_\vph^\flat=P_\vph^{1/2}$.
\eenu
The condition (4) means
the duality between $P_\vph^\al$ and $P_\vph^{1/2-\al}$.
%In particular,
%$P_\vph^\natural$ is self-dual,
%and $(M,H_\vph,P_\vph^\natural,J_\vph)$
%is a standard form.
On the modular involution,
we have
$J_\varphi\xi
=
\Delta_\varphi^{1/2-2\alpha}\xi$
for
$\xi\in P_\varphi^\alpha$.
This shows that
$J_\varphi P_\varphi^\alpha=P_\varphi^{1/2-\alpha}$,
that is,
$J_\vph$ induces
an inversion in the middle point $1/4$.
(See also \cite{miura} for details.)

We set $\mat(\cA_\varphi):=\cA_\varphi\otimes\mat$
and $\varphi_n:=\varphi\otimes\tr_n$.
Then $\mat(\cA_\varphi)$ is a full left Hilbert algebra 
in $\mat(H_\varphi)$.
The multiplication and the involution are given by
\[
[\xi_{i, j}]\cdot[\eta_{i, j}]:=\sum_{k=1}^n[\xi_{i, k}\eta_{k, j}]
\quad\text{and}\quad
[\xi_{i, j}]^\sharp:=[\xi_{j, i}^\sharp]_{i, j}.
\]
Then we have $S_{\varphi_n}=S_\varphi\otimes J_{\tr}$.
Hence the modular operator is 
$\Delta_{\varphi_n}=\Delta_\varphi\otimes\mathrm{id}_{\mat}$.
Denote by $P_{\varphi_n}^{\alpha}$ the positive cone in $\mat(H_\varphi)$
for $\alpha\in[0, 1/2]$. 
We generalize the complete positivity
presented in Definition \ref{defn:cpop}.

\bdf
Let $\alpha\in [0, 1/2]$.
A bounded linear operator $T$ on $H_\varphi$ is said to be 
{\em completely positive} %({\em c.p.})
{\em with respect to $P_\varphi^\alpha$} 
if
$(T\otimes1_{\mat})P_{\varphi_n}^\alpha\subset P_{\varphi_n}^\alpha$
for all $n\in\N$.
\edf

%%%%%%%%%%%%%%%%%%%%%%%%%%%%%%%%%%%%%%%%%%%%%%%%%%%%%%%%%%%%%%%%%%%%%%%%%%%%%%%

\subsection{Completely positive operators from completely positive maps}

%%%%%%%%%%%%%%%%%%%%%%%%%%%%%%%%%%%%%%%%%%%%%%%%%%%%%%%%%%%%%%%%%%%%%%%%%%%%%%%

%We denote by $(H_\varphi, \xi_\varphi)$ 
%the GNS construction of $(M, \varphi)$,
%that is, $\xi_\vph$ is the GNS cyclic vector
%and $H_\vph$ is the norm closure of $M\xi_\vph$.
%We also denote by $\Delta_\varphi$ and $J_\varphi$
%the modular operator and
%the modular conjugation, respectively. 
%The associated modular automorphism group is denoted 
%by $(\sigma_t^\varphi)_{t\in\R}$ as usual.

Let $M$ be a von Neumann algebra
and $\vph\in W(M)$.
Let $C>0$ and $\Phi$ a normal c.p.\ map on $M$
such that 
\begin{equation}
\label{eq:cp}
\varphi\circ\Phi(x)\leq C\varphi(x)
\quad \text{for }
x\in M^+.
\end{equation}
In this subsection,
we will show that $\Phi$ extends
to a c.p.\ operator on $H_\varphi$
with respect to $P_\varphi^\alpha$
for each $\alpha\in[0, 1/2]$.
We use the following folklore among specialists.
(See, for example, \cite[Lemma 4]{ara} for its proof.)

\blem
\label{lem:ara}
Let $T$ be a positive self-adjoint operator 
on a Hilbert space.
For $0\leq r\leq 1$ 
and $\xi\in D(T)$, the domain of $T$, 
we have
$\|T^r\xi\|^2\leq\|\xi\|^2+\|T\xi\|^2$.
\elem

The proof of the following lemma is
inspired by arguments due to Hiai and Tsukada
in \cite[Lemma 2.1]{ht}.

\blem\label{lem:ht}
For $\alpha\in [0, 1/2]$,
one has
\[
\|\Delta_\varphi^\alpha\La_\vph(\Phi(x))\|
\leq C^{1/2}\|\Phi\|^{1/2}
\|\Delta_\varphi^\alpha \La_\vph(x)\|
\quad \text{for }
x\in n_\vph\cap n_\vph^*.
\]
%\[
%\|\Delta_\varphi^\alpha\Phi(x)\xi_\varphi\|
%\leq C^{1/2}\|\Phi\|^{1/2}
%\|\Delta_\varphi^\alpha x\xi_\varphi\|\quad \text{for } 
%x\in M.
%\]
\elem

\bpf
Note that if $x\in n_\vph$,
then $\Ph(x)\in n_\vph$
because
\[
\vph(\Ph(x)^*\Ph(x))
\leq
\|\Ph\|\vph(\Ph(x^*x))
\leq
C\|\Ph\|\vph(x^*x)<\infty
\]

Let $x, y\in n_\vph$ be entire elements
with respect to $\sigma^\varphi$.
We define the entire function $F$ by
\[
F(z):=
\langle
\La_\vph(\Phi(\sigma_{iz/2}^\varphi(x))),
\La_\vph(\sigma_{-i\overline{z}/2}^\varphi(y))
\rangle
\quad
\mbox{for }
z\in \C.
\]
For any $t\in \R$, we have
\begin{align*}
|F(it)|
&=
|\langle
\La_\vph(\Phi(\sigma_{-t/2}^\varphi(x))),
\La_\vph(\sigma_{-t/2}^\varphi(y))
\rangle| \\
&\leq
\|\La_\vph(\Phi(\sigma_{-t/2}^\varphi(x)))\| 
\cdot\|\La_\vph(\sigma_{-t/2}^\varphi(y))\| \\
&=
\vph
(
\Phi(\sigma_{-t/2}^\varphi(x))^*\Phi(\sigma_{-t/2}^\varphi(x))
)^{1/2}
\cdot\|\La_\vph(y)\| \\
%&\leq(
%\|\Phi\|\langle
%\Phi(\sigma_{-t/2}^\varphi(x^*x))\xi_\varphi,
%\xi_\varphi
%\rangle
%)^{1/2}
%\cdot\|y\xi_\varphi\|
%\quad
%\text{by Schwarz inequality}\\
%&\leq(C\|\Phi\|
%\langle
%\sigma_{-t/2}^\varphi(x^*x)\xi_\varphi,
%\xi_\varphi
%\rangle
%)^{1/2}
%\cdot\|y\xi_\varphi\| \quad\text{by (\ref{eq:cp})}
%\\
&\leq C^{1/2}\|\Phi\|^{1/2}\|\La_\vph(x)\|\|\La_\vph(y)\|,
\end{align*}
and
\begin{align*}
|F(1+it)|
&=
|\langle
\Delta_\varphi^{1/2}
\La_\vph(\Phi(\sigma_{(i-t)/2}^\varphi(x))),
\Delta_\varphi^{-it/2}\La_\vph(y)
\rangle| \\
&=
|\langle
J_\varphi
\La_\vph(\Phi(\sigma_{(i-t)/2}^\varphi(x))^*),
\Delta_\varphi^{-it/2}\La_\vph(y)
\rangle| \\
&\leq
\|\La_\vph(\Phi(\sigma_{(i-t)/2}^\varphi(x))^*)\|
\cdot\|\La_\vph(y)\| \\
&=
\vph(
\Phi(\sigma_{(i-t)/2}^\varphi(x))\Phi(\sigma_{(i-t)/2}^\varphi(x))^*)^{1/2}
\cdot\|\La_\vph(y)\| \\
%&\leq(\|\Phi\|
%\langle
%\Phi(\sigma_{(i-t)/2}^\varphi(x)\sigma_{(i-t)/2}^\varphi(x)^*)\xi_\varphi,
%\xi_\varphi
%\rangle
%)^{1/2}
%\cdot\|y\xi_\varphi\| \\
%& \hspace{7cm} \text{by Schwarz inequality}\\
&\leq
C^{1/2}\|\Phi\|^{1/2}
\vph(
\sigma_{(i-t)/2}^\varphi(x)\sigma_{(i-t)/2}^\varphi(x)^*
)^{1/2}
\cdot
\|\La_\vph(y)\| \quad\text{by (\ref{eq:cp})} \\
&=C^{1/2}\|\Phi\|^{1/2}
\|\La_\vph(\sigma_{(i-t)/2}^\varphi(x)^*)\|
\cdot\|\La_\vph(y)\| \\
%&=C^{1/2}\|\Phi\|^{1/2}
%\|\sigma_{\frac{-i-t}{2}}^\varphi(x^*)\xi_\varphi\|
%\cdot\|y\xi_\varphi\| \\
%&=C^{1/2}\|\Phi\|^{1/2}
%\|\Delta_\varphi^{1/2}\sigma_{-\frac{t}{2}}^\varphi(x^*)\xi_\varphi\|
%\cdot\|y\xi_\varphi\| \\
&=C^{1/2}\|\Phi\|^{1/2}
\|J_\varphi\La_\vph(\sigma_{-t/2}^\varphi(x))\|
\cdot\|\La_\vph(y)\| \\
&=C^{1/2}\|\Phi\|^{1/2}
\|\La_\vph(x)\|\|\La_\vph(y)\|.
\end{align*}
Hence the three-lines theorem implies
the following inequality for $0\leq s\leq 1$:
\[
|\langle
\Delta_\varphi^{s/2}
\La_\vph(\Phi(\sigma_{is/2}^\varphi(x))),
\La_\vph(y)
\rangle|
=|F(s)|
\leq C^{1/2}\|\Phi\|^{1/2}
\|\La_\vph(x)\|\|\La_\vph(y)\|.
\]
By replacing $x$ by $\sigma_{-is/2}^\varphi(x)$, we obtain
\[
|\langle
\Delta_\varphi^{s/2}\La_\vph(\Phi(x)),
\La_\vph(y)
\rangle|
\leq C^{1/2}\|\Phi\|^{1/2}
\|\La_\vph(\si_{-is/2}^\vph (x))\|
\|\La_\vph(y)\|.
\]
Since $y$ is an arbitrary entire element of $M$ 
with respect to $\sigma^\varphi$,
we have
\begin{equation}
\label{eq:delta}
\|\Delta_\varphi^{s/2}\La_\vph(\Phi(x))\|
\leq C^{1/2}\|\Phi\|^{1/2}
\|\La_\vph(\sigma_{-is/2}^\varphi(x))\|
=C^{1/2}\|\Phi\|^{1/2}
\|\Delta_\varphi^{s/2}\La_\vph(x)\|.
\end{equation}

For $x\in \cA_\vph$,
take a sequence of entire elements $x_n$ of $M$
with respect to $\sigma^\varphi$
such that
\[
\|\La_\vph(x_n)-\La_\vph(x)\|\to 0
\text{ and }
\|\La_\vph(x_n^*)-\La_\vph(x^*)\|\to 0
\quad
(n\to\infty).
\]
Then we also have
\begin{align*}
\|\Delta_\varphi^{s/2}
\La_\vph(x_n-x)\|^2
&\leq\|\La_\vph(x_n-x)\|^2
+\|\Delta_\varphi^{1/2}\La_\vph(x_n-x)\|^2
\quad\text{by Lemma \ref{lem:ara}} 
\\
&=\|\La_\vph(x_n-x)\|^2+\|\La_\vph(x_n^*-x^*)\|^2 \\
&\to 0.
\end{align*}
Since 
\begin{align*}
\|\La_\vph(\Phi(x_n))-\La_\vph(\Phi(x))\|^2
&=\|\La_\vph(\Phi(x_n-x))\|^2 \\
%&=
%\langle
%\Phi(x_n-x)^*\Phi(x_n-x)\xi_\varphi, 
%\xi_\varphi
%\rangle \\
%&\leq
%\|\Phi\|\langle
%\Phi((x_n-x)^*(x_n-x))\xi_\varphi, 
%\xi_\varphi
%\rangle \quad\text{by Schwarz inequality}\\
%&\leq C\|\Phi\|
%\langle
%(x_n-x)^*(x_n-x)\xi_\varphi, 
%\xi_\varphi
%\rangle \quad\text{by (\ref{eq:cp})} \\
&\leq
C\|\Phi\|
\|\La_\vph(x_n-x)\|^2\to 0,
\end{align*}
we have
\begin{equation}\label{eq:weak}
\langle
\Delta_\varphi^{s/2}\La_\vph(\Phi(x_n)),
\La_\vph(y)
\rangle
\to
\langle
\Delta_\varphi^{s/2}\La_\vph(\Phi(x)),
\La_\vph(y)
\rangle
\quad
\text{for }
y\in n_\vph.
\end{equation}
Moreover, since
\begin{align*}
\|\Delta_\varphi^{s/2}\La_\vph(\Phi(x_m))
-\Delta_\varphi^{s/2}\La_\vph(\Phi(x_n))\|
&\leq
C^{1/2}\|\Phi\|^{1/2}\|\Delta_\varphi^{s/2}\La_\vph(x_m-x_n)\|
\quad\text{by (\ref{eq:delta})}\\
&\to 0\quad (m, n\to\infty),
\end{align*}
the sequence
$\Delta_\varphi^{s/2}\La_\vph(\Phi(x_n))$
is a Cauchy sequence.
Thus
$\Delta_\varphi^{s/2}\La_\vph(\Phi(x_n))$
converges to
$\Delta_\varphi^{s/2}\La_\vph(\Phi(x))$
in norm by (\ref{eq:weak}).
%because $\Delta_\varphi^{\frac{s}{2}}$ is a closed operator.
Therefore, we have
\begin{align*}
\|\Delta_\varphi^{s/2}\La_\vph(\Phi(x))\|
&=
\lim_{n\to\infty}\|\Delta_\varphi^{s/2}\La_\vph(\Phi(x_n))\|
\\
&\leq C^{1/2}\|\Phi\|^{1/2}
\lim_{n\to\infty}\|\Delta_\varphi^{s/2}\La_\vph(x_n)\|
\\
&=C^{1/2}\|\Phi\|^{1/2}
\|\Delta_\varphi^{s/2}\La_\vph(x)\|.
\end{align*}
\epf

%Thanks to Lemma \ref{lem:ht}, the following definition makes sense:

\blem
\label{lem:cp-operator}
Let $M$ be a von Neumann algebra 
with $\vph\in W(M)$ and $\Ph$ be a normal c.p.\ map on $M$.
Suppose $\vph\circ\Ph\leq C\vph$ as before.
Then for $\alpha\in[0, 1/2]$, 
one can define
the bounded operator $T_\Phi^\alpha$ on $H_\varphi$
with $\|T_\Phi^\alpha\|\leq C^{1/2}\|\Phi\|^{1/2}$
by
\[
T_\Phi^\alpha(\Delta_\varphi^\alpha \La_\vph(x))
:=\Delta_\varphi^\alpha \La_\vph(\Phi(x))
\quad\text{for }
x\in n_\vph\cap n_\vph^*.
\]
\elem

It is not hard to see that $T_\Phi^\alpha$ 
in the above
is c.p.\ with respect to $P_\varphi^\alpha$
since $T_\Phi^\al\oti1_{\mat}=T_{\Ph\oti\id_{\mat}}^\al$
preserves $P_{\vph_n}^\al$.
%\bdf
%\label{defn:alHAP}
%Let $M$ be a $\si$-finite von Neumann algebra
%and $\vph$ a faithful normal state on $M$.
%Let $\alpha\in [0, 1/2]$.
%We will say that $M$ has
%the $\al$-\emph{Haagerup approximation property}
%($\al$-HAP)
%if there exists a net of c.c.p.\ operators $T_n$
%on $H$.
%\edf

%%%%%%%%%%%%%%%%%%%%%%%%%%%%%%%%%%%%%%%%%%%%%%%%%%%%%%%%%%%%%%%%%%%%%%%%%%%%%%%%%%%%%%%%%%%%%%%%%%%%%%%%%%%%%%%%%%%%%%%%%%%%%%%%%%%%%%%%%%%%%%%%%%%%%%%%%%%%%%%%%%%%%%%%%%%%%%%%%%%%%%%%%%%%%%%%%%%%%%%%%%%%%%%%%%%%%%%%%%%%%%%%%%%%%%%%%%%%%

\subsection{Haagerup approximation property associated with positive cones}

%%%%%%%%%%%%%%%%%%%%%%%%%%%%%%%%%%%%%%%%%%%%%%%%%%%%%%%%%%%%%%%%%%%%%%%%%%%%%%%%%%%%%%%%%%%%%%%%%%%%%%%%%%%%%%%%%%%%%%%%%%%%%%%%%%%%%%%%%%%%%%%%%%%%%%%%%%%%%%%%%%%%%%%%%%%%%%%%%%%%%%%%%%%%%%%%%%%%%%%%%%%%%%%%%%%%%%%%%%%%%%%%%%%%%%%%%%%%%

We will introduce the ``interpolated'' HAP for a von Neumann algebra.

\bdf\label{def:HAP}
Let $\alpha\in [0, 1/2]$
and $M$ a von Neumann algebra
with $\vph\in W(M)$.
We will say that
$M$ has the $\alpha$-{\em Haagerup approximation property
with respect to $\vph$}
($\al$-HAP$_\vph$) 
if there exists
a net of compact contractive operators $T_n$ on $H_\varphi$
such that $T_n\to 1_{H_{\varphi}}$ in the strong topology
and each $T_n$ is c.p.\ with respect to $P_\varphi^\alpha$.
\edf

%By definition,
%the $1/4$-HAP$_\vph$ is nothing but the HAP
%that is introduced in Definition \ref{defn:HAP}.
%In particular,
%it does not depend on a choice of $\vph$.
%So, we simply say the $1/4$-HAP for $1/4$-HAP$_\vph$.
%Let us state our main theorem in this section.

We will show the above approximation property
is actually a weight-free notion in what follows.

\blem
\label{lem:reduction}
Let $\alpha\in [0, 1/2]$.
Then the following statements hold:
\benu
\item
Let $e\in M_\vph$ be a projection.
If $M$ has the $\al$-HAP$_\vph$,
then $eMe$ has the $\al$-HAP$_{\vph_e}$;

\item
If there exists an increasing net
of projections $e_i$
in $M_\vph$
such that
$e_i\to1$ in the strong topology
and
$e_i Me_i$ has the $\al$-HAP$_{\vph_{e_i}}$
for all $i$,
then $M$ has the $\al$-HAP$_\vph$.
\eenu
\elem

\bpf
(1)
We will regard $H_{\vph_e}=eJeJ H_\vph$,
$J_{\vph_e}=eJe$
and $\De_{\vph_e}=eJeJ\De_\vph$
as usual.
Then it is not so difficult to show that
$P_{\vph_e}^\al=eJeJ P_\vph^\al$.
Take a net $T_n$ as in Definition \ref{def:HAP}.
Then the net $eJeJT_n eJeJ$ does the job.

(2)
Let $\cF$ be a finite subset of $H_\vph$ and $\vep>0$.
Take $i$ such that
\[
\|e_i J_\vph e_i J_\vph\xi-\xi\|<\vep/2
\quad
\mbox{for all }
\xi\in\cF.
\]
We identify $H_{\vph_{e_i}}$ with $e_i J_\vph e_i J_\vph H_\vph$
as usual.
Then take a compact contractive operator
$T$ on $H_{\vph_{e_i}}$
such that
it is c.p.\ with respect to $P_{\vph_{e_i}}^\al$
and satisfies
\[
\|Te_i J_\vph e_i J_\vph\xi-e_i J_\vph e_i J_\vph\xi\|
<\vep/2
\quad
\mbox{for all }
\xi\in \cF.
\]
Thus we have
$\|Te_iJ_\vph e_i J_\vph\xi-\xi\|<\vep$
for $\xi\in \cF$.
It is direct to show that
$Te_iJ_\vph e_i J_\vph$ is a compact contractive operator
such that it is c.p.\ with respect to $P_\vph^\al$,
and we are done.
\epf

\blem
\label{lem:HAPindep}
The approximation property introduced
in Definition \ref{def:HAP}
does not depend on the choice of an f.n.s.\ weight.
Namely, let $M$ be a von Neumann algebra and $\vph,\ps\in W(M)$.
If $M$ has the $\al$-HAP$_\vph$,
then $M$ has the $\al$-HAP$_\ps$.
\elem

\bpf
Similarly as in the proof of Lemma \ref{lem:CS-indep},
it suffices to check that each operation below
inherits the approximation
property introduced in Definition \ref{def:HAP}.
\begin{enumerate}
\item
$\vph\mapsto \vph\oti\Tr$,
where $\Tr$ denotes the canonical tracial weight
on $\bB(\ell^2)$;

\item
$\vph\mapsto \vph_e$,
where $e\in M_\vph$ is a projection;

\item
$\vph\mapsto \vph\circ\al$,
$\al\in \Aut(M)$;

\item
$\vph\mapsto \vph_h$,
where $h$ is a non-singular positive operator
affiliated with $M_\vph$.
\end{enumerate}

(1)
Let $N:=M\oti B(\ell^2)$ and $\ps:=\vph\oti\Tr$.
Take an increasing sequence
of finite rank projections $e_n$ on $\ell^2$
such that $e_n\to1$ in the strong topology.
Then $f_n:=1\oti e_n$ belongs to $N_\ps$
and $f_n Nf_n=M\oti e_n \bB(\ell^2)e_n$,
which has the $\al$-HAP$_{\ps_{f_n}}$.
%because the complete positivity
%of involved operators is assumed.
By Lemma \ref{lem:reduction} (2),
$N$ has the $\al$-HAP$_\ps$.

(2)
This is nothing but Lemma \ref{lem:reduction} (1).

(3).
Let $\ps:=\vph\circ\al$.
Regard as $H_{\ps}=H_\vph$ by putting $\La_\ps=\La_\vph\circ\al$.
We denote by $U_\al$ the canonical unitary implementation, 
which maps $\La_\vph(x)$ to $\La_\ps(\al^{-1}(x))$
for $x\in n_\vph$.
Then it is direct to see that
$\De_\ps=U_\al\De_\vph U_\al^*$,
and $P_\ps^\al=U_\al P_\vph^\al$.
We can show $M$ has the $\al$-HAP$_\ps$
by using $U_\al$.

(4).
Our proof requires a preparation.
We will give a proof after proving Lemma \ref{lem:Th}.
\epf

Let $\al\in[0,1/2]$ and $\vph\in W(M)$.
Note that for an entire element $x\in M$ with respect to $\si^\vph$,
an operator $xJ_\vph \si_{i(\al-\hal)}^\vph(x)J_\vph$ is c.p.\
with respect to $P_\vph^\al$.

\blem
\label{lem:eiejT}
Let $T$ be a c.p.\ operator with respect to $P_\vph^\al$
and $\{e_i\}_{i=1}^m$ a partition of unity in $M_\vph$.
Then the operator
$\sum_{i,j=1}^m e_i J_\vph e_j J_\vph T e_i J_\vph e_j J_\vph$
is c.p.\ with respect to $P_\vph^\al$.
\elem

\begin{proof}
Let $E_{ij}$ be the matrix unit of $\M_m(\C)$.
Set $\rho:=\sum_{i=1}^m e_i\oti E_{1i}$.
Note that $\rho$ belongs to $(M\oti \M_m(\C))_{\vph\oti\tr_n}$.
Then the operator 
\[
\rho J_{\vph\oti\tr_n}\rho J_{\vph\oti\tr_n}
(T\oti1_{\mat})
\rho^* J_{\vph\oti\tr_n}\rho^* J_{\vph\oti\tr_n}
\]
on $H_\vph\oti\M_m(\C)$
is positive with respect to $P_{\vph\oti\tr_n}^\al$
since so is $T\oti1_{\mat}$.
By direct calculation,
this operator equals
$\sum_{i,j=1}^m e_i J_\vph e_j J_\vph T e_i J_\vph e_j J_\vph
\oti E_{11}J_{\tr} E_{11} J_{\tr}$.
Thus we are done.
\end{proof}

Let $h\in M_\vph$ be positive and invertible.
We can put
$\La_{\vph_h}(x):=\La_\vph(xh^{1/2})$
for $x\in n_{\vph_h}=n_\vph$.
This immediately implies that
$\De_{\vph_h}=hJ_\vph h^{-1} J_\vph \De_\vph$,
and
$P_{\vph_h}^\al=h^\al J_\vph h^{\hal} J_\vph P_\vph^\al$.
Thus we have the following result.

\blem
\label{lem:Th}
Let $h\in M_\vph$ be positive and invertible.
If $T$ is a c.p.\ operator with respect to $P_\vph^\al$,
then
\[
T_h:=h^\al J_\vph h^{\hal} J_\vph T h^{-\al} J_\vph h^{-\hal} J_\vph
\]
is c.p.\ with respect to $P_{\vph_h}^\al$.
\elem

\begin{proof}[Resumption of Proof of Lemma \ref{lem:HAPindep}]
Let $\ps:=\vph_h$ and $e(\cdot)$ the spectral resolution of $h$.
Put $e_n:=e([1/n, n])\in M_\vph$ for $n\in\N$.
Note that $M_\vph=M_\ps$.
Since $e_n\to1$ in the strong topology,
it suffices to show that $e_nMe_n$
has the $\al$-HAP$_{\vph_{he_n}}$.
Thus we may and do assume that $h$ is bounded and invertible.

Let us identify $H_\ps=H_\vph$
by putting $\La_\ps(x):=\La_\vph(xh^{1/2})$
for $x\in n_\vph$ as usual,
where we should note that $n_\vph=n_\ps$.
Then we have $\De_\ps=hJ_\vph h^{-1} J_\vph\De_\vph$
and $P_\ps^\al=h^\al J_\vph h^{\hat{\al}}J_\vph P_\vph^\al$ as well.

%Take a net of compact c.c.p.\ operators $T_n$
%with respect to $P_\vph^\al$
%such that $T_n\to1$ strongly.
%Then the net
%$h^\al J_\vph h^{1/2-\al}T_n h^{-\al} J_\vph h^{-(1/2-\al)}$
%does the job.

Let $\cF$ be a finite subset of $H_\vph$ and $\e>0$.
Take $\de>0$ so that $1-(1+\de)^{-1/2}<\e/2$.
Let $\{e_i\}_{i=1}^m$ be a spectral projections of $h$
such that
$\sum_{i=1}^m e_i=1$ and $h e_i \leq \la_i e_i \leq (1+\de)he_i$
for some $\la_i>0$.
Note that $e_i$ belongs to $M_\vph\cap M_{\vph_h}$.
For a c.p.\ operator $T$ with respect to $P_\vph^\al$,
we put
\begin{align*}
T_{h,\de}
&:=
\sum_{i,j=1}^m
e_i J_\vph e_j J_\vph T_h e_i J_\vph e_j J_\vph
\\
&=
\sum_{i,j=1}^m
h^\al e_i J_\vph h^\hal e_j J_\vph T h^{-\al} e_i J_\vph h^{-\hal}e_j J_\vph,
\end{align*}
which is c.p.\ with respect to $P_{\vph_h}^\al$
by Lemma \ref{lem:eiejT} and Lemma \ref{lem:Th}.
The norm of $T_{h,\de}$
equals the maximum of
$\|h^\al e_i Jh^\hal e_j J\,T h^{-\al} e_i Jh^{-\hal}e_j J\|$.
Since we have
\begin{align*}
\|h^\al e_i Jh^\hal e_j J\,T h^{-\al} e_i Jh^{-\hal}e_j J\|
&\leq
\|h^\al e_i\| \|h^\hal e_j\| \|T\| \|h^{-\al} e_i\| \|h^{-\hal}e_j \|
\\
&\leq
\la_i^\al\la_j^\hal
((1+\vep)\la_i^{-1})^{\al}
((1+\vep)\la_j^{-1})^{\hal}
\|T\|
\\
&=
(1+\de)^{1/2},
\end{align*}
we get $\|T_{h,\de}\|\leq (1+\de)^{1/2}$.

Since $M$ has the $\al$-HAP$_\vph$,
we can find a c.c.p.\ compact operator $T$
with respect to $P_\vph^\al$
such that
$\|T_{h,\de}\xi-\xi\|<\e/2$ for all $\xi\in\cF$.
Then $\widetilde{T}:=(1+\de)^{-1/2}T_{h,\de}$ is a c.c.p.\ operator
with respect to $P_{\vph_h}^\al$,
which satisfies
$\|\widetilde{T}\xi-\xi\|<\e$ for all $\xi\in F$.
Thus we are done.
\end{proof}

Therefore,
the $\al$-HAP$_\vph$ does not depend on a choice of $\vph\in W(M)$.
So, we will simply say $\al$-HAP for $\al$-HAP$_\vph$.

Now we are ready to introduce the main theorem in this section.

\bthm\label{thm:equiv}
Let $M$ be a von Neumann algebra.
Then the following statements are equivalent:
\benu
\item
$M$ has the HAP, i.e., the $1/4$-HAP;

\item
$M$ has the $0$-HAP;

\item
$M$ has the $\al$-HAP
for any $\al\in[0,1/2]$;

\item
$M$ has the $\al$-HAP
for some $\al\in[0,1/2]$.

\item
$M$ has the CS-HAP;
\eenu
\ethm

%As an immediate corollary,
%we have the following.
%
%\bcor
%The $\al$-HAP$_\vph$ is a weight-free property.
%Moreover,
%the $\al$-HAP$_\vph$ is equivalent to the HAP.
%\ecor

We will prove the above theorem in several steps.

\begin{proof}[Proof of {\rm (1)$\Rightarrow$(2)} 
in Theorem \ref{thm:equiv}]
Suppose that $M$ has the $1/4$-HAP.
Take an increasing net of $\si$-finite projections $e_i$ in $M$
such that $e_i\to1$ in the strong topology.
Thanks to Lemma \ref{lem:reduction},
it suffices to show that $e_i M e_i$ has the $0$-HAP.
Hence we may and do assume that $M$ is $\si$-finite.
Let $\vph\in M_*^+$ be a faithful state.
By Theorem \ref{thm:sigma-finite},
we can take a net of normal c.c.p.\ maps
$\Ph_n$ on $M$
with
$\vph\circ\Ph_n\leq\vph$
such that 
the following implementing operator $T_n$
is compact and $T_n\to 1_{H_\varphi}$ in the strong topology:
\[
T_n(\De_\vph^{1/4}x\xi_\vph)
=
\De_\vph^{1/4}\Ph_n(x)\xi_\vph
\quad
\mbox{for }
x\in M.
\]
Let $T_{\Ph_n}^0$ be the closure of
$\De_\vph^{-1/4}T_n\De_\vph^{1/4}$
as in Lemma \ref{lem:cp-operator}.
Recall that $T_{\Ph_n}^0$ satisfies
\[
T_{\Ph_n}^0( x\xi_\vph)
=
\Ph_n(x)\xi_\vph
\quad
\mbox{for }
x\in M.
\]
However, the compactness of $T_{\Ph_n}^0$ is not clear.
Thus we will perturb $\Ph_n$ by averaging $\si^\vph$.
Let us put
\[
g_\be(t):=\sqrt{\frac{\be}{\pi}}\exp(-\be t^2)
\quad\text{for }
\be>0
\ \text{and}\ 
t\in\R,
\]
and
\[
U_\be:=\int_\R g_\be(t)\De_\vph^{it}\,dt
=
\widehat{g}_\be(-\log\De_\vph),
\]
where
\[
\widehat{g}_\be(t)
:=
\int_\R g_\be(s)e^{-ist}\,ds
=
\exp(-t^2/(4\be))
\quad
\text{for }
t\in\R.
\]
Then $U_\be\to1$ in the strong topology as $\be\to\infty$.
%By definition, it is obvious that
%$U_\be$ is a c.c.p.\ operator with respect to $P_\vph^{1/4}$.

For $\be,\ga>0$,
we define 
\[
\Phi_{n,\beta, \gamma}(x):=
(\sigma_{g_\beta}^\varphi\circ\Phi_n\circ\sigma_{g_\gamma}^\varphi)(x)
\quad\text{for }
x\in M.
\]
Since $\int_\R g_\gamma(t)\,dt=1$
and $g_\gamma\geq 0$,
the map $\Phi_{n, \beta, \gamma}$ is normal c.c.p.\ 
such that 
$\varphi\circ\Phi_{n, \beta, \gamma}\leq \varphi$.
By Lemma \ref{lem:cp-operator},
we obtain
%for $\alpha\in[0, 1/2]$, 
the associated operator
$T_{\Phi_{n,\beta, \gamma}}^0$,
which is given by
\[
T_{\Phi_{n,\beta, \gamma}}^0
(x\xi_\varphi)
=
\Phi_{n,\beta, \gamma}(x)\xi_\varphi
\quad\text{for }
x\in M.
\]
%Note that
%\[
%\Delta_\varphi^{-\alpha}
%T_\alpha^{\Phi_{n,\beta, \gamma}}\Delta_\varphi^\alpha(x\xi_\varphi)
%=\Phi_{n,\beta, \gamma}(x\xi_\varphi) 
%=T_0^{\Phi_{n,\beta, \gamma}}(x\xi_\varphi)
%\quad\text{for }
%x\in M.
%\]
%Hence
%$\Delta_\varphi^{-\alpha}T_\alpha^{\Phi_{n,\beta, \gamma}}\Delta_\varphi^\alpha$
%extends to the contractive operator on $H_\varphi$,
%which coincides with $T_0^{\Phi_{n,\beta, \gamma}}$.
Moreover, we have
%\begin{equation}
%\label{eq:TU}
$T_{\Phi_{n,\beta, \gamma}}^0
%(x\xi_\varphi)
=
%\Delta_\varphi^{-\alpha}
U_\be
T_{\Ph_n}^0
U_\ga
=
U_\be\De_\vph^{-1/4}
T_n
\De_\vph^{1/4}U_\ga
$.
%\Delta_\varphi^\alpha
%(x\xi_\varphi)
%\quad\text{for }
%x\in M.
%\end{equation}
Hence $T^0_{\Phi_{n,\beta, \gamma}}$ is compact,
because
$e^{-t/4}\widehat{g}_\be(t)$ and 
$e^{t/4}\widehat{g}_\ga(t)$ are bounded functions on $\R$.
Thus we have shown that
$\left(
T_{\Phi_{n,\beta, \gamma}}^0
\right)_{(n,\be,\ga)}$
is a net of contractive compact operators.

It is trivial that
$T^0_{\Phi_{n,\beta, \gamma}}\to1_{H_\vph}$
in the weak topology, 
because $U_\be,U_\ga\to1_{H_\vph}$ as $\be,\ga\to\infty$
and $T_n\to1_{H_\vph}$ as $n\to\infty$ in the strong topology.
\end{proof}

In order to prove Theorem \ref{thm:equiv} (2)$\Rightarrow$(3),
we need a few lemmas.
In what follows,
let $M$ be a von Neumann algebra
with $\vph\in W(M)$.

\blem\label{lem:JTJ}
Let $\alpha\in[0,1/2]$.
Then $M$ has the $\alpha$-HAP$_\vph$
if and only if
$M$ has the $\hat{\al}$-HAP$_\vph$.
\elem

\bpf
It immediately follows from the fact that
$T$ is c.p.\ with respect to $P_\varphi^\alpha$
if and only if 
$J_\varphi TJ_\varphi$ is c.p.\
with respect to $P_\varphi^{\hat{\al}}$.
\epf

\blem\label{lem:uniform}
Let $(U_t)_{t\in \R}$ be a one-parameter unitary group 
and $T$ be a compact operator on a Hilbert space $H$.
If a sequence $(\xi_n)$ in $H$ converges to $0$ weakly,
then $(TU_t\xi_n)$ converges to $0$ in norm,
compact uniformly for $t\in\R$.
\elem

\bpf
Since $T$ is compact,
the map $\R\ni t\mapsto TU_t\in \bB(H)$
is norm continuous.
In particular, for any $R>0$,
the set $\{TU_t\mid t\in[-R,R]\}$ is norm compact.
Since $(\xi_n)$ converges weakly,
it is uniformly norm bounded.
Thus the statement holds
by using a covering of $\{TU_t\mid t\in[-R,R]\}$
by small balls.
\epf

\blem
\label{lem:Pbe}
Let $\al\in[0,1/4]$ and $\be\in[\al,\hat{\al}]$.
Then $P_\vph^\al\subs D(\De_\vph^{\be-\al})$
and $P_\vph^\be=\ovl{\De_\vph^{\be-\al}P_\vph^\al}$.
\elem

\bpf
Since $P_\vph^\al\subs D(\De_\vph^{1/2-2\al})$
and $0\leq \be-\al\leq 1/2-2\al$,
it turns out that
$P_\vph^\al\subs D(\De_\vph^{\be-\al})$.
Let $\xi\in P_\vph^\al$ and take a sequence
$\xi_n\in P_\vph^\sharp$
such that $\De_\vph^\al\xi_n\to\xi$.
Then we have
\begin{align*}
\|\Delta_\varphi^{\be}(\xi_m-\xi_n)\|^2
%\|\Delta_\varphi^{1/4}(x_m-x_n)\xi_\varphi\|^2
&=
\|\Delta_\varphi^{\be-\alpha}
\Delta_\varphi^\alpha (\xi_m-\xi_n)\|^2
%\|\Delta_\varphi^{1/4-\alpha}
%\Delta_\varphi^\alpha (x_m-x_n)\xi_\varphi\|^2
\\
&\leq
\|\De_\vph^0\cdot\Delta_\varphi^\alpha(\xi_m-\xi_n)\|^2
\\
&\quad
+\|\Delta_\varphi^{1/2-2\alpha}\cdot
\Delta_\varphi^\alpha(\xi_m-\xi_n)\|^2 
\quad\text{by Lemma \ref{lem:ara}}
\\
&=
\|\Delta_\varphi^\alpha(\xi_m-\xi_n)\|^2
+\|J_\varphi \Delta_\varphi^\alpha S_\vph(\xi_m-\xi_n)\|^2
\\
&=2\|\Delta_\varphi^\alpha (\xi_m-\xi_n)\|^2
\to 0.
\end{align*}
Hence
$\Delta_\varphi^{\be}\xi_n$ converges to a vector
$\eta$ which belongs to $P_\vph^\be$.
Since
$\Delta_\varphi^{\be-\alpha}(\Delta_\varphi^\alpha \xi_n)
=\Delta_\varphi^{\be}
\xi_n\to \eta$
and
$\Delta_\varphi^{\be-\alpha}$ is closed,
$\Delta_\varphi^{\be-\alpha}\xi=\eta\in P_\varphi^{\be}$.
Hence
$P_\vph^\be\supset \ovl{\De_\vph^{\be-\al}P_\vph^\al}$.
The converse inclusion is obvious
since $\De_\vph^\be P_\vph^\sharp
=\De_\vph^{\be-\al}(\De_\vph^\al P_\vph^\sharp)$.
\epf

%From now on, 
%we assume that $M$ is $\sigma$-finite
%with a faithful state $\varphi\in M_*^+$.
%Let $\vph$ be an f.n.s.\ weight on $M$.

Note that the real subspace 
$R_\varphi^\alpha:=P_\varphi^\alpha-P_\varphi^\alpha$ 
in $H_\varphi$ is closed and the mapping
\[
S_\varphi^\alpha\colon
R_\varphi^\alpha+iR_\varphi^\alpha\ni\xi+i\eta
\mapsto \xi-i\eta\in R_\varphi^\alpha+iR_\varphi^\alpha
\]
is a conjugate-linear closed operator
which has the polar decomposition
\[
S_\varphi^\alpha=J_\varphi\Delta_\varphi^{1/2-2\alpha}.
\]
(See \cite[Poposition 2.4]{ko1} in the case where $M$ is $\sigma$-finite.)

%Next, we shall prove the following result:

\blem
\label{lem:onlyif}
Let $\alpha\in [0, 1/4]$
and $T\in \bB(H_\vph)$ a c.p.\ operator 
with respect to $P_\varphi^\alpha$.
Let $\be\in[\al,\hat{\al}]$.
Then the following statements hold:
\benu
\item
Then the operator
$\Delta_\varphi^{\be-\alpha}T\Delta_\varphi^{\alpha-\be}$
extends to the bounded operator on $H_\varphi$,
which is denoted by $T^\be$ in what follows,
so that $\|T^\be\|\leq\|T\|$.
Also, $T^\be$ is a c.p.\ operator
with respect to $P_\varphi^{\be}$;

\item
If a bounded net of c.p.\ operators $T_n$ 
with respect to $P_\vph^\al$ weakly converges to $1_{H_\vph}$,
then so does the net $T_n^\be$;

\item
If $T$ in {\rm (1)} is non-zero compact,
then so does $T^\be$.
\eenu
\elem

\bpf%[Proof of Lemma \ref{lem:onlyif}]
(1)
Let $\zeta\in P_\vph^{\sharp}$
and $\eta:=\De_\vph^\be \zeta$ which belongs to $P_\vph^\be$.
We put
$\xi:=
T\Delta_\varphi^{\alpha-\be}
\eta
%(\Delta_\varphi^{1/4}x\xi_\varphi).
%=T\Delta_\varphi^\alpha x\xi_\varphi
$.
Since
$\Delta_\varphi^{\alpha-\be}\eta=\De_\vph^\al\zeta\in P_\vph^\al$
and
$T$ is c.p.\
with respect to $P_\varphi^\alpha$,
we obtain $\xi \in P_\varphi^\alpha$.
By Lemma \ref{lem:Pbe},
we know that $\De_\vph^{\be-\al}\xi\in P_\vph^\be$.
Thus
$\Delta_\varphi^{\be-\alpha}T\Delta_\varphi^{\alpha-\be}$
maps $\De_\vph^\be P_\vph^\sharp$
into $P_\varphi^{\be}$.

Hence the complete positivity with respect to $P_\varphi^{\be}$ 
immediately follows
when we prove the norm boundedness of that map.
The proof given below is quite similar as in the one of Lemma \ref{lem:ht}.
Recall the associated Tomita algebra $\mathcal{T}_\vph$.
Let $\xi,\eta\in \mathcal{T}_\vph$.
We define the entire function $F$ by
\[
F(z):=
\langle T \De_\vph^{-z}\xi,
\De_\vph^{\ovl{z}}\eta\rangle
\quad\text{for } z\in\C.
%\langle T\sigma_{iz}^\varphi(x)\xi_\varphi, 
%\sigma_{-i\overline{z}}^\varphi(y)\xi_\varphi\rangle
%\quad\text{for } z\in\C.
\]
For any $t\in\R$, we have
\begin{align*}
|F(it)|
&=
|\langle T\De_\vph^{-it}\xi, 
\De_\vph^{-it}\eta\rangle|
%|\langle T\sigma_{-t}^\varphi(x)\xi_\varphi, 
%\sigma_{-t}^\varphi(y)\xi_\varphi\rangle| \\
\leq 
\|T\|\|\xi\|\|\eta\|.
\end{align*}
Note that 
\begin{align*}
\De_\vph^{-(\hat{\al}-\al+it)}\xi
&=
\Delta_\varphi^\alpha
\De_\vph^{-(\hat{\al}+it)}\xi
\\
&=
\Delta_\varphi^\alpha \xi_1+i\Delta_\varphi^\alpha \xi_2
\in R_\varphi^\alpha+i R_\varphi^\alpha,
\end{align*}
where $\xi_1, \xi_2\in R_\vph^\al$ satisfies
$\De_\vph^{-(\hat{\al}+it)}\xi=\xi_1+i\xi_2$.
Note that $\xi_1$ and $\xi_2$ also belong to $\mathcal{T}_\vph$.
Since $T$ is c.p.\ with respect to $P_\varphi^\alpha$,
we see that $T R_\vph^\al\subs R_\vph^\al$.
Then we have
\begin{align*}
%\Delta_\varphi^{\beta}T\sigma_{i(\beta+it)}^\varphi(x)\xi_\varphi
\Delta_\varphi^{\hat{\al}-\al}T\De_\vph^{-(\hat{\al}-\al+it)}\xi
&=
\Delta_\varphi^{1/2-2\al}
(T\Delta_\varphi^\alpha \xi_1+i T\Delta_\varphi^\alpha \xi_2)
\\
%\Delta_\varphi^{1/2-2\alpha}
%(T\Delta_\varphi^\alpha x_1\xi_\varphi+iT\Delta_\varphi^\alpha x_2\xi_\varphi) \\
&=J_\varphi(T\Delta_\varphi^\alpha \xi_1-iT\Delta_\varphi^\alpha \xi_2) \\
&=J_\varphi T S_\varphi^\alpha
(\Delta_\varphi^\alpha \xi_1+i\Delta_\varphi^\alpha \xi_2) \\
&=J_\varphi T J_\varphi\Delta_\varphi^{1/2-2\alpha}
\De_\vph^{-(\hat{\al}-\al+it)}\xi\\
&=J_\varphi T J_\varphi
\De_\vph^{-it}\xi.
\end{align*}
In particular,
$\Delta_\varphi^{\hat{\al}-\al}T\De_\vph^{-(\hat{\al}-\al)}$ is norm bounded,
and its closure is $J_\varphi T J_\varphi$.
Hence
\begin{align*}
\left|F(\hat{\al}-\al+it)\right|
&=
|\langle T\De_\vph^{-(\hat{\al}-\al+it)}\xi,
\De_\vph^{\hat{\al}-\al-it}\eta\rangle|
\\
&=
|\langle J_\vph T J_\vph \De_\vph^{-it}\xi,\De_\vph^{it}\eta\rangle|
\\
&\leq 
\|T\|\|\xi\|\|\eta\|.
\end{align*}
Applying the three-lines theorem to $F(z)$ at $z=\be-\al\in[0,\hat{\al}-\al]$,
we obtain
\begin{equation}
\label{eq:Tal}
|\langle 
\Delta_\varphi^{\be-\alpha}T\Delta_\varphi^{\alpha-\be}
\xi,\eta\rangle|
=|F(\be-\al)|
\leq\|T\|\|\xi\|\|\eta\|.
\end{equation}
This implies
\[
\|(\Delta_\varphi^{\be-\alpha}T\Delta_\varphi^{\alpha-\be})
 \xi\|
\leq \|T\|\|\xi\|
\quad
\mbox{for all }\xi\in \mathcal{T}_\vph.
\]
Therefore 
$\Delta_\varphi^{\be-\alpha}T\Delta_\varphi^{\alpha-\be}$ 
extends to a bounded operator, which we denote by $T^\be$,
on $H_\varphi$
such that $\|T^\be\|\leq \|T\|$.

(2)
By (1),
we have $\|T_n^\be\|\leq\|T_n\|$,
and thus the net $(T_n^\be)_n$ is also bounded.
Hence the statement follows from the following inequality
for all $\xi,\eta\in \mathcal{T}_\vph$:
\[
|\langle 
(T_n^\be-1_{H_\vph})
\xi,\eta\rangle|
=
|\langle 
(T_n-1_{H_\vph})
\Delta_\varphi^{\alpha-\be}
\xi,\Delta_\varphi^{\be-\al}\eta\rangle|.
\]

(3)
Suppose that $T$ is compact.
Let $\eta_n$ be a sequence in $H_\varphi$ with $\xi_n\to 0$ weakly.
Take $\xi_n\in \mathcal{T}_\vph$
such that $\|\xi_n-\eta_n\|<1/n$ for $n\in\N$.
%Take entire elements $x_n\in M$ with respect to $\sigma^\varphi$
%such that $\|\xi_n-x_n\xi_\varphi\|<1/n$.
It suffices to check that $\|T^\be\xi_n\|\to 0$.
Since the sequence $\xi_n$ is weakly converging, 
there exists $D>0$ such that 
\begin{equation}\label{eq:bound}
\|\xi_n\|\leq D
\quad\text{for all }
n\in\N.
\end{equation}
Let $\eta\in \mathcal{T}_\vph$.
For each $n\in\N$, 
we define the entire function $F_n$ by
\[
F_n(z):=\exp(z^2)
\langle 
T\De_\varphi^{-z}\xi_n,
\De_\vph^{\overline{z}}\eta
\rangle.
\]
Let $\e>0$. Take $t_0>0$ such that 
\begin{equation}\label{eq:exp}
e^{-t^2}\leq\frac{\e}{D\|T\|}
\quad\text{for }
|t|>t_0.
\end{equation}
We let $I:=[-t_0,t_0]$.
Since $T$ is compact,
there exists $n_0\in\N$ such that
\begin{equation}\label{eq:compact}
\|T\Delta_\varphi^{-it}\xi_n\|
\leq \e
\quad\text{and}\quad
\|J_\varphi TJ_\varphi \Delta_\varphi^{-it}\xi_n\|\leq \e
\quad\text{for }
n\geq n_0
\ \text{and}\
t\in I.
\end{equation}
Then for $n\geq n_0$ we have
\begin{align*}
|F_n(it)|
&=e^{-t^2}
|\langle 
T\De_\vph^{-it}\xi_n,
\De_\vph^{-it}\eta
\rangle| \\
&\leq
e^{-t^2}
\|T\Delta_\varphi^{-it}\xi_n\|
\|\eta\|.
\end{align*}
Hence if $t\not\in I$, then
\begin{align*}
|F_n(it)|
&\leq
e^{-t^2}\|T\|
\|\xi_n\|
\|\eta\| \\
&\leq e^{-t^2}D\|T\|
\|\eta\| \quad\text{by}\ (\ref{eq:bound}) \\
&\leq \e\|\eta\|  \quad\text{by}\ (\ref{eq:exp}),
\end{align*}
and if $t\in I$, then
\begin{align*}
|F_n(it)|
&\leq\|T\Delta_\varphi^{-it}\xi_n\|
\|\eta\| \\
&\leq \e\|\eta\| \quad\text{by}\ (\ref{eq:compact}).
\end{align*}
We similarly obtain
\[
|F_n(\hat{\al}-\al+it)|\leq \e\|\eta\|
\quad\text{for }
n\geq n_0
\ \text{and}\
t\in\R.
\]

Therefore the three-lines theorem implies
\[
e^{(\be-\alpha)^2}
|\langle T^\be \xi_n,
\eta\rangle|
=\left|F_n\left(\be-\alpha\right)\right|
\leq\e\|\eta\|
\quad\text{for }
n\geq n_0.
\]
Hence we have
$\|T^\be \xi_n\|\leq\e$
for
$n\geq n_0$.
Therefore $T^\be$ is compact.
\epf

\blem
\label{lem:albe}
Let $M$ be a von Neumann algebra
and $\al\in[0,1/4]$.
If $M$ has the $\al$-HAP,
then
$M$ also has the $\be$-HAP
for all $\be\in[\al,\hat{\al}]$.
\elem
\bpf
Take a net of c.c.p.\ compact operators
$T_n$ with respect to $P_\vph^\al$
as before.
By Lemma \ref{lem:onlyif},
we obtain a net of c.c.p.\ compact operators
$T_n^\be$ with respect to $P_\vph^\be$
such that $T_n^\be$ is converging to $1_{H_\vph}$
in the weak topology.
Thus we are done.
\epf

Now we resume to prove Theorem \ref{thm:equiv}.

\bpf[Proof of {\rm (2)$\Rightarrow$(3)} in Theorem \ref{thm:equiv}]
It follows from Lemma \ref{lem:albe}.
\epf

\bpf[Proof of {\rm (3)$\Rightarrow$(4)} in Theorem \ref{thm:equiv}]
This is a trivial implication.
\epf

\bpf[Proof of {\rm (4)$\Rightarrow$(1)} in Theorem \ref{thm:equiv}]
Suppose that $M$ has the $\alpha$-HAP 
for some $\alpha\in[0, 1/2]$. 
By Lemma \ref{lem:JTJ}, 
we may and do assume that $\alpha\in[0, 1/4]$.
By Lemma \ref{lem:albe},
$M$ has the $1/4$-HAP.
%Employing the previous lemma,
%we see that $M$ has the 0-HAP.
%Then by Lemma \ref{lem:albe},
%$M$ has the $\al$-HAP for all $\al\in[0,1/2]$.
\epf

Therefore we prove the conditions from (1) to (4) are equivalent.
Finally we check the condition (5) and the others are equivalent.

\bpf[Proof of {\rm (1)$\Rightarrow$(5)} in Theorem \ref{thm:equiv}]
It also follows from the proof of {\rm (1)$\Rightarrow$(2)}.
\epf

\bpf[Proof of {\rm (5)$\Rightarrow$(1)} in Theorem \ref{thm:equiv}] 
We may assume that $M$ is $\sigma$-finite
by \cite[Lemma 4.1]{cs} and \cite[Proposition 3.5]{ot}.
Let $\varphi\in M_*^+$ be a faithful state.
For every finite subset $F\subset M$,
we denote by $M_F$ the von Neumann subalgebra
generated by $1$ and
\[
\{\sigma_t^\varphi(x)\mid x\in F, t\in \Q\}.
\]
Then $M_F$ is a separable $\sigma^\varphi$-invariant
and contains $F$.
By \cite[Theorem IX.4.2]{t2}, 
there exists a normal conditional expectation 
$\mathcal{E}_F$ of $M$ onto $M_F$ 
such that $\varphi\circ\mathcal{E}_F=\varphi$.
As in the proof of \cite[Theorem 3.6]{ot}, 
the projection $E_F$ 
on $H_\varphi$ defined below is a c.c.p.\ operator:
\[
E_F(x\xi_\varphi)=\mathcal{E}_F(x)\xi_\varphi
\quad\text{for}\ 
x\in M.
\]
It is easy to see that $M_F$ has the CS-HAP.
It also can be checked that 
if $M_F$ has the HAP for every $F$, 
then $M$ has the HAP.
Hence we can further assume that $M$ is separable.

Since $M$ has the CS-HAP,
there exists a {\em sequence} of normal c.p.\ maps $\Ph_n$
with 
$\vph\circ\Ph_n\leq\vph$
such that
the following implementing operator $T_n^0$ is compact
and $T_n^0\to 1_{H_\vph}$ strongly:
\[
T_n^0(x\xi_\vph):=\Ph_n(x)\xi_\vph
\quad
\mbox{for}\
x\in M.
\]
In particular,
$T_n^0$ is a c.p.\ operator
with respect to $P_\vph^\sharp$.
By the principle of uniform boundedness, 
the sequence $(T_n^0)$ is uniformly norm-bounded.
By Lemma \ref{lem:onlyif},
we have a uniformly norm-bounded sequence of compact operators
$T_n$ such that 
each $T_n$ is c.p.\ with respect to $P_\varphi^{1/4}$
and $T_n$ weakly converges to $1_{H_\vph}$.
By convexity argument,
we may assume that $T_n\to1_{H_\vph}$ strongly.
It turns out from \cite[Theorem 4.9]{ot}
that $M$ has the HAP.
\epf

Therefore we have finished proving Theorem \ref{thm:equiv}.
We will close this section with the following result
that is the contractive map version of Definition \ref{defn:CS1}.

\bthm
\label{thm:cpmap}
Let $M$ be a von Neumann algebra.
Then the following statements are equivalent:
\begin{enumerate}
\item[{\rm (1)}]
$M$ has the HAP;

\item[{\rm (2)}]
For any $\vph\in W(M)$,
there exists a net of normal c.c.p.\ maps
$\Ph_n$ on $M$
such that
\begin{itemize}
\item
$\vph\circ\Ph_n\leq\vph$;

\item
$\Phi_n\to\mathrm{id}_{M}$ in the point-ultraweak topology;

\item
For all $\al\in[0,1/2]$,
the associated c.c.p.\ operators $T_n^\al$ on $H_{\vph}$
defined below are compact
and
$T_n^\al\to\mathrm1_{H_{\vph}}$ in the strong topology:
\begin{equation}
\label{eq:PhTnal}
T_n^\al\Delta_{\vph}^\al \La_\vph(x)
=\Delta_{\vph}^\al\La_\vph(\Phi_n(x))
\quad
\mbox{for all }
x\in n_\vph.
\end{equation}
\end{itemize}

\item[{\rm (3)}]
For some $\vph\in W(M)$ and some $\al\in[0,1/2]$,
there exists a net of normal c.c.p.\ maps
$\Ph_n$ on $M$
such that
\begin{itemize}
\item
$\vph\circ\Ph_n\leq\vph$;

\item
$\Phi_n\to\mathrm{id}_{M}$ in the point-ultraweak topology;

\item
The associated c.c.p.\ operators $T_n^\al$ on $H_{\vph}$
defined below are compact
and
$T_n^\al\to\mathrm1_{H_{\vph}}$ in the strong topology:
\begin{equation}
\label{eq:PhTn}
T_n^\al\Delta_{\vph}^\al \La_\vph(x)
=\Delta_{\vph}^\al\La_\vph(\Phi_n(x))
\quad
\mbox{for all }
x\in n_\vph.
\end{equation}
\end{itemize}
\end{enumerate}
\ethm

First, we will show that
the second statement does not depend on a choice of $\vph$.
So, let us here denote by the approximation property $(\al,\vph)$, 
this approximation property
and by the approximation property $(\al)$ afterwards as well.

\blem
The approximation property $(\al, \vph)$
does not depend on any choice of $\vph\in W(M)$.
\elem

\begin{proof}
Suppose that $M$ has the approximation property $(\al, \vph)$.
It suffices to show that each operation
listed in the proof of Lemma \ref{lem:CS-indep}
inherits the property $(\al,\vph)$.
It is relatively easy to treat the first three operations,
and let us omit proofs for them.
Also,
we can show that
if $e_i$ is a net as in statement of Lemma \ref{lem:reduction} (2)
and $e_i Me_i$ has the approximation property $(\al, \vph_{e_i})$ for each $i$,
then $M$ has the approximation property $(\al, \vph)$.

Thus it suffices to treat $\ps:=\vph_h$ for a positive
invertible element $h\in M_\vph$.
Our idea is similar as in the one of the proof of Lemma \ref{lem:HAPindep}.

Let $\e>0$.
Take $\de>0$ so that $2\de/(1+\de)<\e$.
Let $\{e_i\}_{i=1}^m$ be a spectral projections of $h$
such that
$\sum_{i=1}^m e_i=1$ and $h e_i \leq \la_i e_i \leq (1+\de)he_i$
for some $\la_i>0$.

For a normal c.c.p.\ map $\Ph$ on $M$ such that $\vph\circ\Ph\leq\vph$,
we let
$\Ph_h(x):=h^{-1/2}\Ph(h^{1/2}xh^{1/2})h^{-1/2}$
for $x\in M$.
Then $\Ph_h$ is a normal c.p.\ map satisfying $\ps\circ\Ph_h\leq\ps$.
Next we let
$\Ph_{(h, \de)}(x):=\sum_{i,j=1}^m e_i\Ph_h(e_ixe_j)e_j$
for $x\in M$.
For $x\in M^+$,
we have
\begin{align*}
\ps(\Ph_{(h,\de)}(x))
&=
\sum_{i=1}^m
\ps(e_i\Ph_h(e_ixe_i))
\leq
\sum_{i=1}^m
\ps(\Ph_h(e_ixe_i))
\\
&\leq
\sum_{i=1}^m
\ps(e_ixe_i))
=\ps(x).
\end{align*}
Also, we obtain
\[
\Ph_{(h,\de)}(1)
=
\sum_{i=1}^m e_i\Ph_h(e_i)e_i
=
\sum_{i=1}^m e_ih^{-1/2}\Ph(he_i)h^{-1/2}e_i,
\]
and the norm of
$\Ph_{(h,\de)}(1)$
equals the maximum of
that of $e_ih^{-1/2}\Ph(he_i)h^{-1/2}e_i$.
Since
\begin{align*}
\|e_ih^{-1/2}\Ph(he_i)h^{-1/2}e_i\|
&\leq
\|e_ih^{-1/2}\|^2\|he_i\|
\leq
(1+\de)\la_i^{-1}\cdot\la_i
\\
&=
1+\de,
\end{align*}
we have $\|\Ps_{\de}\|\leq1+\de$.

Now let $\cF$ be a finite subset
in the norm unit ball of $M$
and $\cG$ a finite subset in $M_*$.
Let $\al\in[0,1/2]$.
By the property $(\al, \vph)$,
we can take a normal c.c.p.\ map $\Ph$ on $M$
such that
$\vph\circ\Ph\leq\vph$,
$|\om(\Ph_{(h,\de)}(x)-x)|<\de$
for all $x\in\cF$ and $\om\in\cG$
and the implementing operator $T^\al$ of $\Ph$
with respect to $P_\vph^\al$
is compact.
Put
$\Ps_{(h,\de)}:=(1+\de)^{-1}\Ph_{(h,\de)}$
that is a normal c.c.p.\ map
satisfying $\ps\circ\Ps_{(h,\de)}\leq\ps$.
Then
we have
$|\om(\Ps_{(h,\de)}(x)-x)|<2\de/(1+\de)<\e$
for all $x\in \cF$ and $\om\in \cG$.

By direct computation,
we see that
the implementing operator of $\Ps_{(h,\vep)}$
with respect to $P_\vph^\al$
is equal to the following operator:
\[
\widetilde{T}:=
(1+\de)^{-1}
\sum_{i,j=1}^m
h^\al e_i J_\vph h^\hal e_j J_\vph T h^{-\al} e_i J_\vph h^{-\hal}e_j J_\vph.
\]
Thus $\widetilde{T}$ is compact,
and we are done.
(See also $\widetilde{T}$ in the proof of Lemma \ref{lem:HAPindep}.)
\end{proof}

\begin{proof}[Proof of Theorem \ref{thm:cpmap}]
(1)$\Rightarrow$(2).
%We will show that the HAP implies the strong HAP.
Take $\vph_0\in W(M)$
such that
there exists a partition of unity $\{e_i\}_{i\in I}$
of projections in $M_{\vph_0}$,
the centralizer of $\vph_0$,
such that
$\ps_i:=\vph_0 e_i$ is a faithful normal state on $e_i M e_i$
for each $i\in I$.
Then we have an increasing net of projections $f_j$
in $M_{\vph_0}$ such that $f_j\to1$.
Thus we may and do assume that
$M$ is $\si$-finite as usual.
Employing Theorem \ref{thm:sigma-finite},
we obtain a net of
normal c.c.p.\ maps $\Phi_n$ on $M$ such that
\begin{itemize}
\item
$\varphi\circ\Phi\leq \varphi$;

\item
$\Phi_n\to\mathrm{id}_{M}$ in the point-ultraweak topology;

\item
The operator defined below
is c.c.p.\ compact on $H_\varphi$:
\[
T_n(\Delta_\varphi^{1/4}x\xi_\varphi)
=\Delta_\varphi^{1/4}\Phi_n(x)\xi_\varphi
\ \text{for } x\in M.
\] 
\end{itemize}

Now recall our proof
of Theorem \ref{thm:equiv} (1)$\Rightarrow$(2).
After averaging $\Ph_n$ by $g_\be(t)$ and $g_\ga(t)$,
we obtain a normal c.c.p.\ map $\Ph_{n,\be,\ga}$
which satisfies
$\vph\circ\Ph_{n,\be,\ga}\leq\vph$
and
$\Ph_{n,\be,\ga}\to\id_M$ in the point-ultraweak topology.
For $\al\in[0,1/2]$,
we define the following operator:
\[
T_{n,\be,\ga}^\al\De_\vph^\al\La_\vph(x)
:=
\De_\vph^\al\La_\vph(\Ph_{n,\be,\ga}(x))
\quad
\mbox{for }
x\in n_\vph.
\]
Then we can show the compactness of $T_{n,\be,\ga}^\al$
in a similar way to the proof
of Theorem \ref{thm:equiv} (1)$\Rightarrow$(2),
and we are done.

(2)$\Rightarrow$(3).
This implication is trivial.

(3)$\Rightarrow$(1).
By our assumption,
we have a net of c.c.p.\ compact operators $T_n^\al$
with respect to some $P_\vph^\al$
such that $T_n^\al\to1$ in the strong operator topology.
Namely $M$ has the $\al$-HAP,
and thus $M$ has the HAP by Theorem \ref{thm:equiv}.
\end{proof}

%\bcor[cf.\ {\cite[Theorem 2.3]{cost}}, {\cite[Remark 5.8]{ot}}]
%Let $M$ be a von Neumann algebra. 
%Then $M$ has the HAP in the sense of \cite{ot}
%if and only if
%$M$ has the HAP in the sense of \cite{cs}.
%\ecor

%\bpf
%By considering an increasing net 
%of $\sigma$-finite projections 
%converging to the identity in $M$, 
%we may assume that $M$ is $\sigma$-finite.
%Then the corollary follows from Theorem \ref{thm:equiv}.
%\epf

%\brem
%In \cite{ot}, we introduce 
%the notion of the Haagerup approximation property
%for an arbitrary von Neumann algebra 
%in terms of the standard form.
%In \cite{cs}, Caspers and Skalski 
%independently  give the one
%in terms of c.p.\ maps 
%approximating to the identity
%with respect to a fixed faithful normal semifinite weight.
%After these works,
%two approaches are equivalent 
%by using the permanence property for cores 
%(see \cite{cost} \cite[Remark 5.8]{ot} for details).
%However this proof is quite indirect.
%Now we expect that the presented proof is straightforward. 
%\erem

%%%%%%%%%%%%%%%%%%%%%%%%%%%%%%%%%%%%%%%%%%%%%%%%%%%%%%%%%%%%%%%%%%%%%%%%%%%%%%%%%%%%%%%%%%%%%%%%%%%%%%%%%%%%%%%%%%%%%%%%%%%%%%%%%%%%%%%%%%%%%%%%%%%%%%%%%%%%%%%%%%%%%%%%%%%%%%%%%%%%%%%%%%%%%%%%%%%%%%%%%%%%%%%%%%%%%%%%%%%%%%%%%%%%%%%%%%%%

\section{Haagerup approximation property and non-commutative $L^p$-spaces}

%%%%%%%%%%%%%%%%%%%%%%%%%%%%%%%%%%%%%%%%%%%%%%%%%%%%%%%%%%%%%%%%%%%%%%%%%%%%%%%%%%%%%%%%%%%%%%%%%%%%%%%%%%%%%%%%%%%%%%%%%%%%%%%%%%%%%%%%%%%%%%%%%%%%%%%%%%%%%%%%%%%%%%%%%%%%%%%%%%%%%%%%%%%%%%%%%%%%%%%%%%%%%%%%%%%%%%%%%%%%%%%%%%%%%%%%%%%%

In this section,
we study some relations between the Haagerup approximation property 
and non-commutative $L^p$-spaces 
associated with a von Neumann algebra.

%%%%%%%%%%%%%%%%%%%%%%%%%%%%%%%%%%%%%%%%%%%%%%%%%%%%%%%%%%%%%%%%%%%%%%%%%%%%%%%%%%%%%%%%%%%%%%%%%%%%%%%%%%%%%%%%%%%%%%%%%%%%%%%%%%%%%%%%%%%%%%%%%%%%%%%%%%%%%%%%%%%%%%%%%%%%%%%%%%%%%%%%%%%%%%%%%%%%%%%%%%%%%%%%%%%%%%%%%%%%%%%%%%%%%%%%%%%%

\subsection{Haagerup's $L^p$-spaces}

%%%%%%%%%%%%%%%%%%%%%%%%%%%%%%%%%%%%%%%%%%%%%%%%%%%%%%%%%%%%%%%%%%%%%%%%%%%%%%%%%%%%%%%%%%%%%%%%%%%%%%%%%%%%%%%%%%%%%%%%%%%%%%%%%%%%%%%%%%%%%%%%%%%%%%%%%%%%%%%%%%%%%%%%%%%%%%%%%%%%%%%%%%%%%%%%%%%%%%%%%%%%%%%%%%%%%%%%%%%%%%%%%%%%%%%%%%%%

We begin with Haagerup's $L^p$-spaces in \cite{haa2}. 
(See also \cite{te1}.)
Throughout this subsection,
we fix an f.n.s.\ weight $\varphi$
on a von Neumann algebra $M$.
We denote by $R$ 
the crossed product $M\rtimes_\sigma\R$ of $M$
by the $\R$-action $\sigma:=\sigma^\varphi$.
Via the natural embedding, we have the inclusion $M\subs R$.
Then $R$ admits the canonical faithful normal semifinite trace $\tau$
and there exists the dual action $\theta$ satisfying
$\tau\circ\theta_s=e^{-s}\tau$
for $s\in\R$.
Note that $M$ is equal to the fixed point algebra $R^\th$,
that is,
$M=\{y\in R \mid \theta_s(y)=y\ \text{for } s\in\R\}$.

We denote by $\widetilde{R}$
the set of all $\tau$-measurable
closed densely defined operators
affiliated with $R$.
The set of positive elements in $\widetilde{R}$
is denoted by $\widetilde{R}^+$.
For $\psi\in M_*^+$,
we denote by $\hat{\psi}$ its dual weight
on $R$
and by $h_\psi$ the element of $\widetilde{R}^+$
satisfying
$
\hat{\psi}(y)=\tau(h_\psi y)
$
%\quad\text{for all }
for all $y\in R$.

Then the map $\psi\mapsto h_\psi$ is extended to
a linear bijection of $M_*$
onto the subspace
\[
\{h\in \widetilde{R} \mid \theta_s(h)=e^{-s}h
\ \text{for } s\in\R\}.
\]
Let $1\leq p< \infty$. 
The $L^p$-space of $M$ due to Haagerup is defined as follows:
\[
L^p(M):=
\{
a\in\widetilde{R} 
\mid 
\theta_s(a)=e^{-\frac{s}{p}}a
\ \text{for } 
s\in \R
\}.
\]
Note that the spaces $L^p(M)$ and their relations are
independent of the choice of an f.n.s. weight $\varphi$,
and thus canonically associated with a von Neumann algebra $M$.
Denote by $L^p(M)^+$ the cone $L^p(M)\cap\widetilde{R}^+$.
Recall that if $a\in\widetilde{R}$ with 
the polar decomposition $a=u|a|$,
then $a\in L^p(M)$ 
if and only if 
$|a|^p\in L^1(M)$.
The linear functional $\mathrm{tr}$ on $L^1(M)$ 
is defined by
\[
\mathrm{tr}(h_\psi):=\psi(1)
\quad\text{for }
\psi\in M_*.
\]
Then $L^p(M)$ becomes a Banach space 
with the norm
\[
\|a\|_p:=\mathrm{tr}(|a|^p)^{1/p}
\quad\text{for }
a\in L^p(M).
\]
In particular, $M_*\simeq L^1(M)$
via the isometry $\psi\mapsto h_\psi$.
For non-commutative $L^p$-spaces, 
the usual H\"older inequality also holds.
Namely, let $q>1$ with $1/p+1/q=1$,
and we have
\[
|\tr(ab)|\leq\|ab\|_1\leq\|a\|_p\|b\|_q
\quad\text{for }
a\in L^p(M), b\in L^q(M).
\]
Thus the form $(a, b)\mapsto\tr(ab)$ gives a duality
between $L^p(M)$ and $L^q(M)$. 
Moreover the functional $\tr$ has the ``tracial'' property:
\[
\tr(ab)=\tr(ba)
\quad\text{for }
a\in L^p(M), b\in L^q(M).
\] 
Among non-commutative $L^p$-spaces,
$L^2(M)$ becomes a Hilbert space 
with the inner product
\[
\langle a, b\rangle:=\mathrm{tr}(b^*a)
\quad\text{for }
a, b\in L^2(M).
\]

The Banach space $L^p(M)$ has the natural
$M$-$M$-bimodule structure
as defined below:
\[
x\cdot a\cdot y:=xay
\quad
\mbox{for }
x, y\in M,\ a\in L^p(M).
\]
The conjugate-linear isometric involution $J_p$ on $L^p(M)$
is defined by $a\mapsto a^*$ for $a\in L^p(M)$.
Then the quadruple $(M, L^2(M), J_2, L^2(M)^+)$ 
is a standard form.

%%%%%%%%%%%%%%%%%%%%%%%%%%%%%%%%%%%%%%%%%%%%%%%%%%%%%%%%%%%%%%%%%%%%%%%%%%%%%%%%%%%%%%%%%%%%%%%%%%%%%%%%%%%%%%%%%%%%%%%%%%%%%%%%%%%%%%%%%%%%%%%%%%%%%%%%%%%%%%%%%%%%%%%%%%%%%%%%%%%%%%%%%%%%%%%%%%%%%%%%%%%%%%%%%%%%%%%%%%%%%%%%%%%%%%%%%%%%

\subsection{Haagerup approximation property for non-commutative $L^p$-spaces}

%%%%%%%%%%%%%%%%%%%%%%%%%%%%%%%%%%%%%%%%%%%%%%%%%%%%%%%%%%%%%%%%%%%%%%%%%%%%%%%%%%%%%%%%%%%%%%%%%%%%%%%%%%%%%%%%%%%%%%%%%%%%%%%%%%%%%%%%%%%%%%%%%%%%%%%%%%%%%%%%%%%%%%%%%%%%%%%%%%%%%%%%%%%%%%%%%%%%%%%%%%%%%%%%%%%%%%%%%%%%%%%%%%%%%%%%%%%%

We consider the f.n.s.\ weight 
$\varphi^{(n)}:=\varphi\otimes\mathrm{tr}_n$
on $\mat(M):=M\otimes\mat$.
Since 
$\sigma_t^{(n)}:=\sigma_t^{\varphi^{(n)}}=\sigma_t\otimes\mathrm{id}_n$,
we have 
\[
R^{(n)}:=\mat(M)\rtimes_{\sigma^{(n)}}\R
=(M\rtimes_\sigma\R)\otimes\mat=\mat(R).
\]
The canonical f.n.s.\ trace on $R^{(n)}$ is given by
$\tau^{(n)}=\tau\otimes\mathrm{tr}_n$. 
Moreover $\theta^{(n)}:=\theta\otimes\mathrm{id}_n$ is the dual action 
on $R^{(n)}$.
Since $\widetilde{R^{(n)}}=\mat(\widetilde{R})$, we have
\[
L^p(\mat(M))=\mat(L^p(M))
\quad\text{and}\quad
\mathrm{tr}^{(n)}=\mathrm{tr}\otimes\mathrm{tr}_n.
\]

\bdf
Let $M$ and $N$ be two von Neumann algebras 
with f.n.s.\ weights $\varphi$ and $\psi$, respectively.
For $1\leq p\leq \infty$,
a bounded linear operator $T\colon L^p(M)\to L^p(N)$ 
is {\em completely positive} %({\em c.p.}) 
if $T^{(n)}\colon L^p(\mat(M))\to L^p(\mat(N))$ 
is positive for every $n\in\N$,
where 
$T^{(n)}[a_{i, j}]=[Ta_{i, j}]$
%\quad\text{for }
for $[a_{i, j}]\in L^p(\mat(M))=\mat(L^p(M))$.
\edf

In the case where $M$ is $\sigma$-finite,
the following result gives a construction 
of a c.p.\ operator on $L^p(M)$
from a c.p.\ map on $M$.

\bthm[cf. {\cite[Theorem 5.1]{hjx}}]
\label{thm:hjx}
If $\Phi$ is a c.c.p.\ map on $M$ 
with $\varphi\circ\Phi\leq C\varphi$,
then one obtain a c.p.\ operator $T^p_\Phi$
on $L^p(M)$ 
with $\|T^p_\Phi\|\leq C^{1/p}\|\Phi\|^{1-1/p}$, 
which is defined by
\begin{equation}\label{eq:L^p-op}
T^p_\Phi(h_\varphi^{1/2p}xh_\varphi^{1/2p})
:=h_\varphi^{1/2p}\Phi(x)h_\varphi^{1/2p}
\quad\text{for }
x\in M.
\end{equation}
\ethm

Let $M$ be a $\sigma$-finite von Neumann algebra
with a faithful state $\varphi\in M_*^+$.
Since
\[
\|h_\varphi^{1/4} x h_\varphi^{1/4}\|_2^2 \\
=\mathrm{tr}(h_\varphi^{1/4}x^*h_\varphi^{1/2}xh_\varphi^{1/4}) =\|\Delta_\varphi^{1/4}x\xi_\varphi\|^2
\quad\text{for } x\in M,
\]
we have the isometric isomorphism 
$L^2(M)\simeq H_\varphi$
defined by
$h_\varphi^{1/4}xh_\varphi^{1/4}
\mapsto
\Delta_\varphi^{1/4}x\xi_\varphi$
for $x\in M$.
Therefore under this identification,
the above operator $T_\Phi^2$ is nothing but $T_\Phi^{1/4}$,
which is given in Lemma \ref{lem:cp-operator}.

\bdf
Let $1<p<\infty$
and $M$ be a von Neumann algebra.
We will say that $M$ has the $L^p$-\emph{Haagerup approximation property}
($L^p$-HAP)
if
there exists a net of c.c.p.\ compact operators $T_n$ on $L^p(M)$
such that $T_n\to1_{L^p(M)}$ in the strong topology.
\edf

Note that a von Neumann algebra $M$ has the HAP
if and only if $M$ has the $L^2$-HAP,
because $(M, L^2(M), J_2, L^2(M)^+)$ is a standard form
as mentioned previously.

%%%%%%%%%%%%%%%%%%%%%%%%%%%%%%%%%%%%%%%%%%%%%%%%%%%%%%%%%%%%%%%%%%%%%%%%%%%%%%%%%%%%%%%%%%%%%%%%%%%%%%%%%%%%%%%%%%%%%%%%%%%%%%%%%%%%%%%%%%%%%%%%%%%%%%%%%%%%%%%%%%%%%%%%%%%%%%%%%%%%%%%%%%%%%%%%%%%%%%%%%%%%%%%%%%%%%%%%%%%%%%%%%%%%%%%%%%%%

\subsection{Kosaki's $L^p$-spaces}

%%%%%%%%%%%%%%%%%%%%%%%%%%%%%%%%%%%%%%%%%%%%%%%%%%%%%%%%%%%%%%%%%%%%%%%%%%%%%%%%%%%%%%%%%%%%%%%%%%%%%%%%%%%%%%%%%%%%%%%%%%%%%%%%%%%%%%%%%%%%%%%%%%%%%%%%%%%%%%%%%%%%%%%%%%%%%%%%%%%%%%%%%%%%%%%%%%%%%%%%%%%%%%%%%%%%%%%%%%%%%%%%%%%%%%%%%%%%

We assume that 
$\varphi$ is a faithful normal state 
on a $\sigma$-finite von Neumann algebra $M$.
For each $\eta\in[0,1]$,
$M$ is embedded into $L^1(M)$ by 
$M\ni x\mapsto h_\varphi^\eta x h_\varphi^{1-\eta}\in L^1(M)$.
We define the norm 
$\|h_\varphi^\eta x h_\varphi^{1-\eta}\|_{\infty, \eta}
:=\|x\|$
on $h_\varphi^\eta Mh_\varphi^{1-\eta}\subset L^1(M)$,
i.e., $M\simeq h_\varphi^\eta Mh_\varphi^{1-\eta}$. 
Then $(h_\varphi^\eta Mh_\varphi^{1-\eta}, L^1(M))$ becomes
a pair of compatible Banach spaces
in the sense of A. P. Calder\'{o}n \cite{ca}.
For $1<p<\infty$, 
Kosaki's $L^p$-space $L^p(M; \varphi)_\eta$
is defined as the complex interpolation space
$C_\theta(h_\varphi^\eta Mh_\varphi^{1-\eta}, L^1(M))$
equipped with the complex interpolation norm
$\|\cdot\|_{p, \eta}:=\|\cdot\|_{C_\theta}$,
where $\theta=1/p$.
In particular, $L^p(M; \varphi)_0$, 
$L^p(M; \varphi)_1$ and $L^p(M; \varphi)_{1/2}$
are called the left, the right 
and the symmetric $L^p$-spaces, respectively. 
Note that the symmetric $L^p$-space $L^p(M; \varphi)_{1/2}$
is exactly the $L^p$-space studied in \cite{te2}.

From now on, we assume that $\eta=1/2$,
and we will use the notation $L^p(M; \varphi)$ 
for the symmetric $L^p$-space $L^p(M; \varphi)_{1/2}$.

Note that $L^p(M; \varphi)$ is exactly 
$h_\varphi^{1/2q}L^p(M)h_\varphi^{1/2q}$,
where $1/p+1/q=1$,
and
\[
\|h_\varphi^{1/2q}ah_\varphi^{1/2q}\|_{p, 1/2}=\|a\|_p
\quad\text{for }
a\in L^p(M).
\]
Namely, we have 
$L^p(M; \varphi)
=h_\varphi^{1/2q}L^p(M)h_\varphi^{1/2q}
\simeq L^p(M)$.
Furthermore, we have
\[
h_\varphi^{1/2}Mh_\varphi^{1/2}
\subset L^p(M; \varphi)\subset L^1(M),
\]
and $h_\varphi^{1/2}Mh_\varphi^{1/2}$ 
is dense in $L^p(M; \varphi)$.

Let $\Phi$ be a c.p.\ map on $M$ 
with $\varphi\circ\Phi\leq \varphi$.
Note that $T^2_\Phi$ in Theorem \ref{thm:hjx}
equals $T^{1/4}_\Phi$ in Lemma \ref{lem:cp-operator}
under the identification $L^2(M; \varphi)=H_\varphi$.
By the reiteration theorem 
for the complex interpolation method
in \cite{bl,ca},
we have
\begin{equation}\label{eq:p>2}
L^p(M; \varphi)
=C_{2/p}(h_\varphi^{1/2} Mh_\varphi^{1/2}, L^2(M; \varphi))
\quad\text{for } 2<p<\infty,
\end{equation}
and
\begin{equation}\label{eq:p<2}
L^p(M; \varphi)
=C_{\frac{2}{p}-1}(L^2(M; \varphi), L^1(M))
\quad\text{for } 1<p<2.
\end{equation}
(See also \cite[Section 4]{ko3}.)
Thanks to \cite{lp}, 
if $T^2_\Phi=T^{1/4}_\Phi$ is compact 
on $L^2(M; \varphi)=H_\varphi$,
then $T^p_\Phi$ is also compact 
on $L^p(M; \varphi)$ 
for $1<p<\infty$.

%%%%%%%%%%%%%%%%%%%%%%%%%%%%%%%%%%%%%%%%%%%%%%%%%%%%%%%%%%%%%%%%%%%%%%%%%%%%%%%%%%%%%%%%%%%%%%%%%%%%%%%%%%%%%%%%%%%%%%%%%%%%%%%%%%%%%%%%%%%%%%%%%%%%%%%%%%%%%%%%%%%%%%%%%%%%%%%%%%%%%%%%%%%%%%%%%%%%%%%%%%%%%%%%%%%%%%%%%%%%%%%%%%%%%%%%%%%%

\subsection{The equivalence between the HAP and the $L^p$-HAP}

%%%%%%%%%%%%%%%%%%%%%%%%%%%%%%%%%%%%%%%%%%%%%%%%%%%%%%%%%%%%%%%%%%%%%%%%%%%%%%%%%%%%%%%%%%%%%%%%%%%%%%%%%%%%%%%%%%%%%%%%%%%%%%%%%%%%%%%%%%%%%%%%%%%%%%%%%%%%%%%%%%%%%%%%%%%%%%%%%%%%%%%%%%%%%%%%%%%%%%%%%%%%%%%%%%%%%%%%%%%%%%%%%%%%%%%%%%%%

We first show that the HAP implies the $L^p$-HAP
in the case where $M$ is $\sigma$-finite.

\bthm\label{thm:L^p-HAP}
Let $M$ be a $\sigma$-finite von Neumann algebra
with a faithful state $\varphi\in M_*^+$.
Suppose that $M$ has the HAP, i.e.,
there exists a net of normal c.c.p.\ map $\Phi_n$ on $M$
with $\varphi\circ\Phi_n\leq\varphi$
satisfying the following: 
\begin{itemize}
\item $\Phi_n\to\mathrm{id}_M$ in the point-ultraweak topology; 
\item the associated operators $T^2_{\Phi_n}$ on $L^2(M)$
defined below are compact and
$T^2_{\Phi_n}\to 1_{L^2(M)}$ in the strong topology:
\[
T^2_{\Phi_n}(h_\varphi^{1/4}xh_\varphi^{1/4})
=h_\varphi^{1/4}\Phi_n(x)h_\varphi^{1/4}
\quad\text{for }
x\in M.
\]
\end{itemize}
Then $T^p_{\Phi_n}\to 1_{L^p(M)}$ in the strong topology on $L^p(M)$
for $1<p<\infty$.
In particular, $M$ has the $L^p$-HAP for all $1<p<\infty$.
\ethm

\bpf
We will freely use notations and results in \cite{ko3}.
First we consider the case where $p>2$. 
By (\ref{eq:p>2}) we have
\[
L^p(M; \varphi)
=C_\theta(h_\varphi^{1/2} Mh_\varphi^{1/2}, L^2(M; \varphi))
\quad\text{with}\
\theta:=2/p.
\] 
Let $a\in L^p(M; \varphi)$ 
with $\|a\|_{L^p(M; \varphi)}=\|a\|_{C_\theta}\leq 1$ 
and $0<\e<1$. 
By the definition of the interpolation norm,
there exists $f\in F(h_\varphi^{1/2} Mh_\varphi^{1/2}, L^2(M; \varphi))$ 
such that $a=f(\theta)$ 
and $|\!|\!| f |\!|\!|_F\leq 1+\e/3$.
By \cite[Lemma 4.2.3]{bl} (or \cite[Lemma 1.3]{ko3}), 
there exists $g\in F_0(h_\varphi^{1/2} Mh_\varphi^{1/2}, L^2(M; \varphi))$ 
such that $|\!|\!| f-g |\!|\!|_F\leq \e/3$
and $g(z)$ is of the form
\[
g(z)=\exp(\lambda z^2)\sum_{k=1}^K\exp(\lambda_kz)
h_\varphi^{1/2}x_kh_\varphi^{1/2},
\]
where $\lambda>0$, $K\in\N$, 
$\lambda_1,\dots,\lambda_K\in\R$
and $x_1,\dots,x_K\in M$.
Then
\[
\|f(\theta)-g(\theta)\|_\theta\leq|\!|\!| f-g |\!|\!|_F\leq \e/3.
\]

Since 
\[
\lim_{t\to\pm\infty}\|g(1+it)\|_{L^2(M; \varphi)}=0,
\]
a subset $\{g(1+it) \mid t\in\R\}$ of $L^2(M; \varphi)$
is compact in norm.
Hence there exists $n_0\in\N$ such that
\[
\|T^2_{\Phi_n} g(1+it)-g(1+it)\|_{L^2(M; \varphi)}
\leq\left(\frac{\e}{4^{1-\theta}3}\right)^{1/\theta}
\quad\text{for }
n\geq n_0
\ \text{and}\
t\in\R.
\]
Moreover,
\begin{align*}
\|\Phi_n(g(it))-g(it)\|
&\leq\|\Phi_n-\mathrm{id}_M\|\|g(it)\| \\
&\leq 2|\!|\!| g |\!|\!|_F \\
&\leq 2\left(|\!|\!| f |\!|\!|_F+\e/3\right) \\
&\leq 2\left(1+2\vep/3\right)<4.
\end{align*}
We put 
\[
T_{\Phi_n} g(z)
:=\exp(\lambda z^2)\sum_{k=1}^K\exp(\lambda_kz)
h_\varphi^{1/2}\Phi_n(x_k)h_\varphi^{1/2}
\in F_0(h_\varphi^{1/2}Mh_\varphi^{1/2}, L^2(M; \varphi)).
\]
Then $T^p_{\Phi_n} g(\theta)=T_{\Phi_n} g(\theta)\in L^p(M; \varphi)$.
Hence by \cite[Lemma 4.3.2]{bl} (or \cite[Lemma A.1]{ko3}),
we have
\begin{align*}
\|(T^p_{\Phi_n} g)(\theta)-g(\theta)\|_\theta
&\leq
\left(
\int_\R\|\Phi_n(g(it))-g(it)\|P_0(\theta, t)
\,
\frac{dt}{1-\theta}
\right)^{1-\theta} \\
&\quad\times\left(
\int_\R\|T^2_{\Phi_n} g(1+it)-g(1+it)\|_{L^2(M; \varphi)}P_1(\theta, t)
\,\frac{dt}{\theta}
\right)^\theta \\
&\leq 4^{1-\theta}\cdot \e/(4^{1-\theta}3)
=\e/3.
\end{align*}
Therefore since $T^p_{\Phi_n}$ are contractive on $L^p(M; \varphi)$, we have
\begin{align*}
\|T^p_{\Phi_n} f(\theta)-f(\theta)\|_\theta
&\leq \|T^p_{\Phi_n} f(\theta)-T^p_{\Phi_n} g(\theta)\|_\theta
+\|T^p_{\Phi_n} g(\theta)-g(\theta)\|_\theta
\\
&\quad
+\|g(\theta)-f(\theta)\|_\theta \\
&<\e.
\end{align*}
Hence $T^p_{\Phi_n}\to 1_{L^p(M; \varphi)}$ in the strong topology.

In the case where $1<p<2$, 
the same argument also works.
\epf

We continue further investigation of the $L^p$-HAP.

\blem\label{lem:dual}
Let $1<p, q<\infty$ with $1/p+1/q=1$.
Then $M$ has the $L^p$-HAP
if and only if
$M$ has the $L^q$-HAP.
\elem

\bpf
Suppose that $M$ has the $L^p$-HAP, i.e.,
there exists a net of c.c.p.\ compact operators $T_n$
on $L^p(M)$ such that 
$T_n\to1_{L^p(M)}$ in the strong topology.
Then we consider the transpose operators ${}^tT_n$ on $L^q(M)$,
which are defined by
\[
\mathrm{tr}({}^tT_n(b)a)=\mathrm{tr}(bT_n(a))
\quad\text{for }
a\in L^p(M), b\in L^q(M).
\]
It is easy to check that 
${}^tT_n$ is c.c.p.\ compact 
and ${}^tT_n\to 1_{L^q(M)}$ in the weak topology.
By taking suitable convex combinations, we have a net of c.c.p.\ compact operators 
$\widetilde{T}_n$
on $L^q(M)$ such that 
$\widetilde{T}_n\to1_{L^q(M)}$ in the strong topology.
Hence $M$ has the $L^q$-HAP.
\epf

We use the following folklore among specialists.
(See \cite[Proposition 7.6]{pt}, \cite[Poposition 3.1]{ko2}.)

\blem\label{lem:pt}
Let $h$ and $k$ be a $\tau$-measurable self-adjoint operators
such that $h$ is non-singular.
Then there exists $x\in M^+$ 
such that $k=h^{1/2}xh^{1/2}$
if and only if
$k\leq ch$ for some $c\geq 0$.
In this case, we have $\|x\|\leq c$.
\elem

%\bpf
%If $K=H^{1/2}TH^{1/2}$, then $K\leq \|T\|H$. 
%Conversely, if $K\leq c H$, 
%then $K^{1/2}=c^{1/2}SH^{1/2}$
%for some $S\in\bB(\cH)$ with $\|S\|\leq 1$.
%Put $T:=c S^*S$. 
%Then $0\leq T\leq c 1_H$.
%Moreover
%\[
%H^{1/2}TH^{1/2}=cH^{1/2}S^*SH^{1/2}=K.
%\]
%\epf

In the case where $p=2$, 
the following lemma is proved in \cite[Lemma 4.1]{ot}.

\blem
\label{lem:order}
Let $1<p<\infty$
and $M$ be a $\sigma$-finite von Neumann algebra
with $h_0\in L^1(M)^+$ such that 
$h_0^{1/2}$ is cyclic and separating in $L^2(M)$.
Then $\Theta_{h_0}^p\colon M_{\mathrm{sa}}\to L^p(M)$, 
which is defined by
\[
\Theta_{h_0}^p(x):=h_0^{1/2p}xh_0^{1/2p}
\quad\text{for } 
x\in M_{\mathrm{sa}},
\]
induces an order isomorphism between
$\{x\in M_{\mathrm{sa}}\mid -c1\leq x\leq c1\}$
and
$K_{h_0}^p:=\{h\in L^p(M)_{\mathrm{sa}}\mid 
-ch_0^{1/p}\leq a\leq ch_0^{1/p}\}$
for each $c>0$.
Moreover $\Theta_{\xi_0}^p$
is $\sigma(M, M_*)$-$\sigma(L^p(M), L^q(M))$ continuous.
\elem

\bpf
%In the case where $p=2$, the lemma is nothing but \cite[Lemma 4.1]{ot}.
Suppose that $p>2$
and take $q>1$ with $1/p+1/q=1$.
First we will show that $\Theta_{h_0}^p$ is
$\sigma(M, M_*)$-$\sigma(L^p(M), L^q(M))$ continuous. 
If $x_n\to 0$ in $\sigma(M, M_*)$, 
then for $b\in L^q(M)$ we have
\[
\tr(\Theta_{h_0}^p(x_n)b)
=\tr((h_0^{1/2p}x_nh_0^{1/2p})b)
=\tr(x_n(h_0^{1/2p}bh_0^{1/2p}))\to 0,
\]
because $h_0^{1/2p}bh_0^{1/2p}\in L^1(M)=M_*$.

Next we will prove that $\Theta_{h_0}^p$ is an order isomorphism 
between $\{x\in M\mid 0\leq x\leq 1\}$ 
and $\{a\in L^p(M)\mid 0\leq a\leq h_0^{1/p}\}$. 
If $x\in M$ with $0\leq x\leq 1$, then 
\[
\tr((h_0^{1/p}-h_0^{1/2p}xh_0^{1/2p})b) \\
=\tr((1-x)h_0^{1/2p}bh_0^{1/2p})\geq 0
\quad\text{for }
b\in L^q(M)^+.
\]
Hence $h_0^{1/p}\geq \Theta_{h_0}^p(x)=h_0^{1/2p}xh_0^{1/2p}\geq 0$.

Conversely, let $a\in L^p(M)$ with $0\leq a\leq h_0^{1/p}$. 
By Lemma \ref{lem:pt}, 
there exists $x\in M$ with $0\leq x\leq 1$
such that $a=h_0^{1/2p}xh_0^{1/2p}$.
\epf

We will use the following results. 

\blem[{\cite[Theorem 4.2]{ko4}}]
\label{thm:norm-conti}
For $1\leq p, q<\infty$, the map
\[
L^p(M)^+\ni a\mapsto a^{\frac{p}{q}}\in L^q(M)^+
\]
is a homeomorphism with respect to the norm topologies.
\elem

In \cite{ko5}, it was proved that Furuta's inequality \cite{fu} remains valid
for $\tau$-measurable operators. 
In particular, the L\"{o}wner--Heinz inequality holds 
for $\tau$-measurable operators. 

\blem
\label{lem:LH-ineq}
If $\tau$-measurable positive self-adjoint operators $a$ and $b$ 
satisfy $a\leq b$,
then $a^r\leq b^r$ for $0<r<1$.
\elem

The following lemma can be proved 
similarly as in the proof of \cite[Lemma 4.2]{ot}.

\blem
\label{lem:cp-ex}
Let $1\leq p<\infty$. 
If $a\in L^p(M)^+$, then
\benu
\item A functional $f_a\colon L^q(M)\to\C$, $b\mapsto\mathrm{tr}(ba)$
is a c.p.\ operator;
\item An operator $g_a\colon\C\to L^p(M)$, $z\mapsto z a$
is a c.p.\ operator.
\eenu
\elem

In the case where $p=2$,
the following lemma is also proved in \cite[Lemma 4.3]{ot}.
We give a proof for reader's convenience

\blem
\label{lem:cyc-sep}
Let $1< p<\infty$ and 
$M$ be a $\sigma$-finite von Neumann algebra
with a faithful state $\varphi\in M_*^+$.
If $M$ has the $L^p$-HAP,
then there exists a net of c.c.p.\ compact operators $T_n$ 
on $L^p(M)$ 
such that $T_n\to1_{L^p(M)}$ in the strong topology,
and $(T_nh_\varphi^{1/p})^{p/2}\in L^2(M)^+$ 
is cyclic and separating for all $n$.
\elem

\bpf
Since $M$ has the $L^p$-HAP,
there exists a net of c.c.p.\ compact operators $T_n$
on $L^p(M)$ such that 
$T_n\to1_{L^p(M)}$ in the strong topology.
Set $a_n^{1/p}:=T_nh_\varphi^{1/p}\in L^p(M)^+$.
Then $a_n\in L^1(M)^+$. If we set
\[
h_n:=a_n+(a_n-h_\varphi)_-\in L^1(M)^+,
\]
then $h_n\geq h_\varphi$.
By Lemma \ref{lem:LH-ineq},
we obtain $h_n^{1/2}\geq h_\varphi^{1/2}$.
It follows from \cite[Lemma 4.3]{co2} that 
$h_n^{1/2}\in L^2(M)^+$ is cyclic and separating.
Now we define a compact operator $T_n'$ on $L^p(M)$ by
\[
T_n'a:=T_na+\mathrm{tr}(ah_\varphi^{1/q})
(h_n^{1/p}-a_n^{1/p})
\quad\text{for }
a\in L^p(M).
\] 
Since $h_n^{1/p}\geq a_n^{1/p}$ 
by Lemma \ref{lem:LH-ineq},
each $T_n$ is a c.p.\ operator, because of Lemma \ref{lem:cp-ex}.
Note that 
\[
T_n'h_\varphi^{1/p}=
T_nh_\varphi^{1/p}+\mathrm{tr}(h_\varphi)
(h_n^{1/p}-a_n^{1/p})
=h_n^{1/p}.
\]
Since 
$a_n^{1/p}=T_nh_\varphi^{1/p}\to h_\varphi^{1/p}$
in norm, 
we have $a_n\to h_\varphi$ in norm by Lemma \ref{thm:norm-conti}.
Since 
\[
\|h_n-a_n\|_1
=\|(a_n-h_\varphi)_-\|_1
\leq\|a_n-h_\varphi\|_1
\to 0,
\]
we obtain $\|h_n^{1/p}-a_n^{1/p}\|_p\to 0$ 
by Lemma \ref{thm:norm-conti}.
Therefore $\|T_n'a-a\|_p\to 0$ for any $a\in L^p(M)$.
Since $\|T_n'-T_n\|\leq\|h_n^{1/p}- a_n^{1/p}\|_p\to 0$,
we get $\|T_n'\|\to 1$.
Then operators $\widetilde{T}_n:=\|T_n'\|^{-1}T_n'$ give a desired net.
\epf

If $M$ is $\sigma$-finite 
and the $L^p$-HAP 
for some $1<p<\infty$,
then we can recover a net of normal c.c.p.\ maps on $M$
approximating to the identity
such that the associated implementing operators
on $L^p(M)$ are compact.
In the case where $p=2$,
this is nothing but \cite[Theorem 4.8]{ot}
(or Theorem \ref{thm:cpmap}).

\bthm
\label{thm:L^p-HAP2}
Let $1<p<\infty$ and $M$ a $\sigma$-finite von Neumann algebra 
with a faithful state $\varphi\in M_*^+$.
If $M$ has the $L^p$-HAP, then 
there exists a net of normal c.c.p.\ map $\Phi_n$ on $M$
with $\varphi\circ\Phi_n\leq\varphi$
satisfying the following:
\begin{itemize}
\item $\Phi_n\to\mathrm{id}_M$ in the point-ultraweak topology;
\item the associated c.c.p.\ operator $T^p_{\Phi_n}$ on $L^p(M)$
defined below are compact
and $T^p_{\Phi_n}\to1_{L^p(M)}$ in the strong topology:
\[
T^p_{\Phi_n}(h_\varphi^{1/2p}xh_\varphi^{1/2p})
=h_\varphi^{1/2p}\Phi_n(x)h_\varphi^{1/2p}
\quad\text{for }
x\in M.
\]
\end{itemize}
\ethm

\bpf
The case where $p=2$ is nothing but \cite[Theorem 4.8]{ot}.
Let $p\ne 2$.
Take $q>1$ such that $1/p+1/q=1$.
By Lemma \ref{lem:cyc-sep},
there exists a net of c.c.p.\ compact operator $T_n$ 
on $L^p(M)$ 
such that $T_n\to1_{L^p(M)}$ in the strong topology,
and $h_n^{1/2}:=(T_nh_\varphi^{1/p})^{p/2}$ is cyclic and separating 
on $L^2(M)$ for all $n$.

Let $\Theta_{h_\varphi}^p$ and $\Theta_{h_n}^p$ be
the maps given in Lemma \ref{lem:order}.
For each $x\in M_{\mathrm{sa}}$, 
take $c>0$ such that $-c1\leq x\leq c1$.
Then 
\[
-ch_\varphi^{1/p}
\leq h_\varphi^{1/2p}x h_\varphi^{1/2p}
\leq ch_\varphi^{1/p}.
\]
Since $T_n$ is positive,
we have 
\[
-ch_\varphi^{1/p}
\leq T_n(h_\varphi^{1/2p}x h_\varphi^{1/2p})
\leq ch_\varphi^{1/p}.
\]
From Lemma \ref{lem:order},
the operator $(\Theta_{h_n}^p)^{-1}(T_n(h_\varphi^{1/2p}x h_\varphi^{1/2p}))$ in $M$ is well-defined.
Hence we can define a linear map $\Phi_n$ on $M$ by
\[
\Phi_n:=(\Theta_{h_n}^p)^{-1}\circ T_n\circ \Theta_{h_\varphi}^p.
\]
In other words,
\[
T_n(h_\varphi^{1/2p}x h_\varphi^{1/2p})
=h_n^{1/2p}\Phi_n(x) h_n^{1/2p}
\quad\text{for }
x\in M.
\]
One can easily check that $\Phi_n$ is a normal u.c.p.\ map.

\vspace{5pt}
\noindent
\textbf{ Step 1.} We will show that $\Phi_n\to\mathrm{id}_M$ 
in the point-ultraweak topology.
\vspace{5pt}

Since $h_\varphi^{1/2}Mh_\varphi^{1/2}$ 
is dense in $L^1(M)$,
it suffices to show that 
\[
\mathrm{tr}(\Phi_n(x)h_\varphi^{1/2}yh_\varphi^{1/2})\to
\mathrm{tr}(xh_\varphi^{1/2}yh_\varphi^{1/2})
\quad\text{for }
x, y\in M.
\]
However
\begin{align*}
|\mathrm{tr}((\Phi_n(x)-x)h_\varphi^{1/2}yh_\varphi^{1/2})|
&=|\mathrm{tr}(h_\varphi^{1/2p}(\Phi_n(x)-x)h_\varphi^{1/2p}\cdot h_\varphi^{\frac{1}{2q}}yh_\varphi^{\frac{1}{2q}})| \\
&=|\mathrm{tr}((T_n-1_{L_p(M)})(h_\varphi^{1/2p}xh_\varphi^{1/2p})\cdot h_\varphi^{\frac{1}{2q}}yh_\varphi^{\frac{1}{2q}})| \\
&\leq\|(T_n-1_{L_p(M)})(h_\varphi^{1/2p}xh_\varphi^{1/2p})\|_p\|h_\varphi^{\frac{1}{2q}}yh_\varphi^{\frac{1}{2q}}\|_q \\
&\to 0.
\end{align*}

\vspace{5pt}
\noindent
\textbf{ Step 2.} We will make a small perturbation of $\Phi_n$.
\vspace{5pt}

By Lemma \ref{thm:norm-conti}, 
we have $\|h_n-h_\varphi\|_1\to 0$, i.e., 
$\|\varphi_n-\varphi\|\to 0$, 
where $\varphi_n\in M_*^+$ is the unique element 
with $h_n=h_{\varphi_n}$.
By a similar argument as in the proof of \cite[Theorem 4.8]{ot},
one can obtain normal c.c.p.\ maps $\widetilde{\Phi}_n$ on $M$
with $\widetilde{\Phi}_n\to\mathrm{id}_M$ in the point-ultraweak topology,
and c.c.p.\ compact operators $\widetilde{T}_n$ on $L^p(M)$
with $\widetilde{T}_n\to 1_{L^p(M)}$ in the strong topology
such that $\varphi\circ\widetilde{\Phi}_n\leq \varphi$ and 
\[
\widetilde{T}_n(h_\varphi^{1/2p}xh_\varphi^{1/2p})
=h_\varphi^{1/2p}\widetilde{\Phi}_n(x)h_\varphi^{1/2p}
\quad\text{for }
x\in M.
\]
Moreover operators $\widetilde{T}_n$ are nothing but $T^p_{\Phi_n}$.
\epf

Recall that $M$ has the completely positive approximation property (CPAP)
if and only if $L^p(M)$ has the CPAP 
for some/all $1\leq p<\infty$.
This result is proved in \cite[Theorem 3.2]{jrx}.
The following is the HAP version of this fact.

\bthm\label{thm:L^p-HAP3}
Let $M$ be a von Neumann algebra.
Then the following are equivalent:
\benu
\item $M$ has the HAP;
\item $M$ has the $L^p$-HAP for all $1<p<\infty$;
\item $M$ has the $L^p$-HAP for some $1<p<\infty$.
\eenu
\ethm

\bpf
We first reduce the case where $M$ is $\sigma$-finite 
by the following elementary fact
similarly as in the proof of \cite[Theorem 3.2]{jrx}.
Take an f.n.s.\ weight $\varphi$ on $M$ 
and an increasing net of projection $e_n$ in $M$ 
with $e_n\to 1_M$ in the strong topology
such that $\sigma_t^\varphi(e_n)=e_n$ for all $t\in\R$
and $e_nMe_n$ is $\sigma$-finite for all $n$.
Then we can identify $L^p(e_nMe_n)$ with a subspace of $L^p(M)$
and there exists a completely positive projection 
from $L^p(M)$ onto $L^p(e_nMe_n)$ via $a\mapsto e_nae_n$.
Moreover the union of these subspaces is norm dense in $L^p(M)$.
Therefore it suffices to prove the theorem in the case where $M$ is $\sigma$-finite.

(1)$\Rightarrow$(2). 
It is nothing but Theorem \ref{thm:L^p-HAP}.

(2)$\Rightarrow$(3). 
It is trivial.

(3)$\Rightarrow$(1). 
Suppose that $M$ has the $L^p$-HAP for some $1<p<\infty$.
We may and do assume that $p<2$ by Lemma \ref{lem:dual}.
Let $\varphi\in M_*$ be a faithful state.
By Theorem \ref{thm:L^p-HAP2}, 
there exists a net of normal c.c.p.\ maps $\Phi_n$ on $M$
with $\varphi\circ\Phi_n\leq \varphi$
such that $\Phi_n\to\mathrm{id}_M$ in the point-ultraweak topology
and a net of the associated compact operators $T^p_{\Phi_n}$ 
converges to $1_{L^p(M)}$ in the strong topology.
By the reiteration theorem for the complex interpolation method,
we have
$L^2(M; \varphi)
=C_\theta(h_\varphi^{1/2}Mh_\varphi^{1/2}, L^p(M; \varphi))$
for some $0<\theta<1$.
By \cite{lp}, operators $T^2_{\Phi_n}$ are also compact.
Moreover, by the same argument 
as in the proof of Theorem \ref{thm:L^p-HAP},
we have $T^2_{\Phi_n}\to 1_{L^2(M)}$ in the strong topology.
\epf

\brem
In the proof of \cite[Theorem 3.2]{jrx}, 
it is shown that c.p.\ operators on $L^p(M)$ give c.p.\ maps on $M$ 
by using the result of L. M. Schmitt in \cite{sch}. 
See \cite[Theorem 3.1]{jrx} for more details.
However our approach is much different 
and based on the technique of A. M. Torpe in \cite{tor}. 
\erem

%%%%%%%%%%%%%%%%%%%%%%%%%%%%%%%%%%%%%%%%%%%%%%%%%%%%%%%%%%%%%%%%%%%%%%%%%%%%%%%%%%%%%%%%%%%%%%%%%%%%%%%%%%%%%%%%%%%%%%%%%%%%%%%%%%%%%%%%%%%%%%%%%%%%%%%%%%%%%%%%%%%%%%%%%%%%%%%%%%%%%%%%%%%%%%%%%%%%%%%%%%%%%%%%%%%%%%%%%%%%%%%%%%%%%%%%%%%%%


\begin{thebibliography}{DCFY}

\bibitem[Ar]{ara}
H. Araki; 
\textit{Some properties of modular conjugation operator of von Neumann algebras and a non-commutative Radon--Nikodym theorem with a chain rule.} 
Pacific J. Math. \textbf{50} (1974), 309--354. 

\bibitem[AM]{am}
H. Araki and T. Masuda;
\textit{Positive cones and $L^p$-spaces for von Neumann algebras}, 
Publ. Res. Inst. Math. Sci. \textbf{18} (1982), 339--411.


\bibitem[AH]{ah1}
H. Ando and U. Haagerup; 
\textit{Ultraproducts of von Neumann algebras}.
Preprint. arXiv:1212.5457.

\bibitem[BL]{bl}
J. Bergh and J. L\"{o}fstr\"{o}m
\textit{Interpolation spaces. An introduction}. 
Grundlehren der Mathematischen Wissenschaften, No. \textbf{223}. Springer-Verlag, Berlin-New York, 1976. x+207 pp. 

\bibitem[Br1]{br1}
M. Brannan;
\textit{Approximation properties for free orthogonal and free unitary quantum groups.} 
J. Reine Angew. Math. \textbf{672} (2012), 223--251. 

\bibitem[Br2]{br2}
M. Brannan; 
\textit{Reduced operator algebras of trace-preserving quantum automorphism groups.}
Doc. Math. \textbf{18} (2013), 1349--1402.

\bibitem[Ca]{ca}
A.-P. Calder\'{o}n; 
\textit{Intermediate spaces and interpolation, the complex method}. 
Studia Math. \textbf{24} (1964) 113--190. 

\bibitem[CS]{cs}
M. Caspers and A. Skalski;
\textit{The Haagerup property for arbitrary von Neumann algebras}.
Preprint. arXiv:1312.1491.

\bibitem[C+]{cost}
M. Caspers, R. Okayasu, A. Skalski and R. Tomatsu;
\textit{Generalisations of the Haagerup property 
to arbitrary von Neumann algebras}.
To appear in C. R. Acad. Sci. Paris Ser. I Math.

%\bibitem[C+]{book}
%P-A. Cherix, M. Cowling, P. Jolissaint, P. Julg and A. Valette, Alain; 
%{\it Groups with the Haagerup property.} 
%Gromov's a-T-menability. Progress in Mathematics, \textbf{ 197}. Birkh\"{a}user Verlag, Basel, 2001. viii+126 pp.

\bibitem[Ch]{cho}
M. Choda; 
\textit{Group factors of the Haagerup type.} 
Proc. Japan Acad. Ser. A Math. Sci. \textbf{59} (1983), no. 5, 174--177.

\bibitem[Co1]{co1}
A. Connes;
\textit{Une classification des facteurs de type III.}
Ann. Sci. \'{E}cole Norm. Sup. (4) \textbf{6} (1973), 133--252.

\bibitem[Co2]{co2}
A. Connes; 
\textit{Caract\'{e}risation des espaces vectoriels ordonn\'{e}s sous-jacents aux alg\`{e}bres de von Neumann.} 
Ann. Inst. Fourier (Grenoble) \textbf{24} (1974), no. 4, x, 121--155 (1975).

%\bibitem[Co2]{co2}
%A. Connes;
%{\it Classification of injective factors. Cases II$_1$, II$_\infty$, III$_\la$,
%$\la\neq1$.}
%Ann. of Math. (2)  \textbf{ 104} (1976), no. 1, 73--115. 


%\bibitem[CJ]{cj}
%A. Connes and V. Jones; 
%{\it Property T for von Neumann Algebras.}
%Bull. London Math. Soc.  \textbf{ 17}  (1985),  no. 1, 57--62.

\bibitem[CK]{lp}
M. Cwikel and N. J. Kalton;
\textit{Interpolation of compact operators by the methods of Calder\'{o}n and Gastavsson-Peetre.}
Proc. Edinburgh Math. Soc. \textbf{38} (1995), 261--276

\bibitem[D+]{dfsw}
M. Daws, P. Fima, A. Skalski and S. White;
\textit{The Haagerup property for locally compact quantum groups.}
To appear in J. Reine Angew. Math.

\bibitem[DCFY]{cfy}
K. De Commer, A. Freslon and M. Yamashita; 
\textit{CCAP for universal discrete quantum groups.}
To appear in Comm. Math. Phys.

\bibitem[EL]{el}
E. G. Effros and E. C. Lance;
\textit{Tensor products of operator algebras.} 
Adv. Math. \textbf{25} (1977), no. 1, 1--34. 

\bibitem[Fu]{fu}
T. Furuta;
\textit{$A\geq B\geq 0$ assures $(B^rA^p B^r)^{1/q}\geq B^{(p+2r)/q}$
for $r\geq0$, $p\geq 0$, $q\geq 1$
with $(1+2r)q\geq p+2r$}.
Proc. Amer. Math. Soc. \textbf{101} (1987), no. 1, 85--88.

\bibitem[Ha1]{haa1}
U. Haagerup; 
\textit{The standard form of von Neumann algebras.}
Math. Scand. \textbf{37} (1975), no. 2, 271--283.

\bibitem[Ha2]{haa2}
U. Haagerup;
\textit{$L^p$-spaces associated with an arbitrary von Neumann algebra}. 
Alg\`{e}bres d'op\'{e}rateurs et leurs applications en physique math\'{e}matique (Proc. Colloq., Marseille, 1977), pp. 175--184.

\bibitem[Ha3]{haa3}
U. Haagerup; 
\textit{An example of a nonnuclear $C^*$-algebra, which has the metric approximation property.} 
Invent. Math. \textbf{50} (1978/79), no. 3, 279--293.

%\bibitem[Ha3]{haa3}
%U. Haagerup;
%{\it Operator-valued weights in von Neumann algebras. II}.
%J. Funct. Anal. \textbf{ 33} (1979), no. 3, 339--361.


\bibitem[HJX]{hjx}
U. Haagerup, M. Junge and Q. Xu;
\textit{A reduction method for noncommutative $L_p$-spaces and applications}.
Trans. Amer. Math. Soc. \textbf{362} (2010), no. 4, 2125--2165. 

\bibitem[Han]{han}
F. Hansen;
\textit{Les espaces $L^p$ d'une alg\'{e}bre de von Neumann}, 
J. Funct. Anal. \textbf{40} (1981), 151--169.

\bibitem[HT]{ht}
F. Hiai and M. Tsukada;
\textit{Generalized conditional expectations and martingales in noncommutative $L^p$-spaces}. 
J. Operator Theory \textbf{18} (1987), no. 2, 265--288. 

\bibitem[Izu]{izu}
H. Izumi;
\textit{Constructions of non-commutative $L^p$-spaces with a complex parameter arising from modular actions}, 
Internat. J. Math. \textbf{8} (1997), no. 8, 1029--1066. 

\bibitem[Jo]{jol}
P. Jolissaint; 
\textit{Haagerup approximation property for finite von Neumann algebras.} 
J. Operator Theory \textbf{48} (2002), no. 3, suppl., 549--571.

\bibitem[JR]{jr}
M. Junge and Z-J. Ruan;
\textit{Approximation properties for noncommutative $L_p$-spaces associated with discrete groups.} 
Duke Math. J. \textbf{117} (2003), no. 2, 313--341. 

\bibitem[JRX]{jrx}
M. Junge, Z-J. Ruan and Q. Xu;
\textit{Rigid $\mathcal{OL}_p$ structures of non-commutative $L_p$-spaces associated with hyperfinite von Neumann algebras.} 
Math. Scand. \textbf{96} (2005), no. 1, 63--95. 

\bibitem[JX]{jx}
M. Junge and Q. Xu;
\textit{Noncommutative Burkholder/Rosenthal inequalities}. 
Ann. Probab. \textbf{31} (2003), no. 2, 948--995. 

\bibitem[Ko1]{ko1}
H. Kosaki;
\textit{Positive cones associated with a von Neumann algebra}. 
Math. Scand. \textbf{47} (1980), no. 2, 295--307. 

\bibitem[Ko2]{ko2}
H. Kosaki;
\textit{Positive cones and $L^p$-spaces associated with a von Neumann algebra}. 
J. Operator Theory \textbf{6} (1981), no. 1, 13--23. 

\bibitem[Ko3]{ko3}
H. Kosaki;
\textit{Applications of the complex interpolation method to a von Neumann algebra: noncommutative $L^p$-spaces}. 
J. Funct. Anal. \textbf{56} (1984), no. 1, 29--78.

\bibitem[Ko4]{ko4}
H. Kosaki;
\textit{Applications of uniform convexity of noncommutative $L^p$-spaces}.
Trans. Amer. Math. Soc. \textbf{283} (1984), no. 1, 265--282.

\bibitem[Ko5]{ko5}
H. Kosaki;
\textit{A remark on Sakai's quadratic Radon-Nikod\'{y}m theorem}.
Proc. Amer. Math. Soc. \textbf{116} (1992), no. 3, 783--786.

\bibitem[KV]{kv}
J. Kustermans and S. Vaes;
\textit{Locally compact quantum groups in the von Neumann algebraic setting.}
Math. Scand.  \textbf{92} (2003), no. 1, 68--92.

\bibitem[Le]{le}
F. Lemeux;
\textit{Haagerup property for quantum reflection groups.}
To appear Proc. Amer. Math. Soc.

\bibitem[Miu]{miura}
Y. Miura;
\textit{Some properties of the convex cones in a Hilbert space.}
Artes liberales, \textbf{25}, (1979), 165--176.

\bibitem[MT]{mt}
Y. Miura and J. Tomiyama; 
\textit{On a characterization of the tensor product of self-dual cones associated to the standard von Neumann algebras.} 
Sci. Rep. Niigata Univ. Ser. A No. \textbf{20} (1984), 1--11.

\bibitem[OT]{ot}
R. Okayasu and R. Tomatsu;
\textit{Haagerup approximation property for arbitrary von Neumann algebras}.
Preprint. arXiv:1312.1033

\bibitem[PT]{pt}
G. K. Pedersen and M. Takesaki;
\textit{The operator equation $THT=K$}. 
Proc. Amer. Math. Soc. \textbf{36} (1972), 311--312. 

\bibitem[Sch]{sch}
L. M. Schmitt;
\textit{Facial structures on $L^p$-spaces of $W^*$-algebras}, 
Dissertation, Universit\"{a}t des Saarlandes, Saarbr\"{u}cken, 1985.

\bibitem[SW]{sw1}
L. M. Schmitt and G. Wittstock;
\textit{Characterization of matrix-ordered standard forms of $\mathrm{W}^*$-algebras.} 
Math. Scand. \textbf{51} (1982), no. 2, 241--260 (1983).

\bibitem[St]{st}
\c{S}. Str\v{a}til\v{a};
\textit{Modular theory in operator algebras.}
Editura Academiei Republicii Socialiste Rom\^{a}nia,
Bucharest; Abacus Press, Tunbridge Wells, 1981. 492pp.

\bibitem[Ta]{t2}
M. Takesaki; 
\textit{Theory of operator algebras. II.}
Encyclopaedia of Mathematical Sciences, \textbf{ 125}.
Operator Algebras and Non-commutative Geometry, \textbf{6}.
Springer-Verlag, Berlin, 2003. xxii+518 pp.

%\bibitem[Ta3]{t3}
%M. Takesaki; {\it Theory of operator algebras. III.} 
%Encyclopaedia of Mathematical Sciences, \textbf{ 127}. Operator Algebras and Non-commutative Geometry, \textbf{ 8}. Springer-Verlag, Berlin, 2003. xxii+548 pp.

\bibitem[Te1]{te1}
M. Terp;
\textit{$L^p$-spaces associated with von Neumann algebras}.
Notes, K\o benhavns
Universitets Matematiske Institut, 1981.

\bibitem[Te2]{te2}
M. Terp;
\textit{Interpolation spaces between a von Neumann algebra and its predual}. 
J. Operator Theory \textbf{8} (1982), no. 2, 327--360. 

\bibitem[To]{tor}
A. M. Torpe; 
\textit{A characterization of semidiscrete von Neumann algebras in terms of matrix ordered Hilbert spaces.} 
Preprint, Matematisk Institut, Odense Universitet, 1981.

\end{thebibliography}
\end{document}